\documentclass[12pt]{article}
\usepackage{amsfonts}
\usepackage{amsmath}
\usepackage{authblk}
\usepackage{hyperref}

\normalbaselines

\title{A Mean Field Game System and a Related Deterministic Optimal Control Problem}
\author{\c{S}tefana-Lucia Ani\c{t}a} 
\date{\small Octav Mayer Institute of Mathematics of the Romanian Academy, \\ Bd. Carol I 8, Ia\c{s}i 700505, Romania, e-mail: stefi$\_$anita@yahoo.com}

\begin{document}

\maketitle

\begin{abstract} This paper concerns a Mean Field Game (MFG) system related to a Nash type equilibrium for dynamical games associated to large populations.
One shows that the MFG system may be viewed as the Euler-Lagrange system for an optimal control problem related to a Fokker-Planck equation with control in the drift.
One derives the existence of a weak solution to the MFG system and under more restrictive assumptions one proves a uniqueness result. 
\end{abstract}

\section{Introduction}

The Mean Field Game (MFG) theory has been introduced in two seminal works by J.M. Lasry and P.L. Lions \cite{LL} for the mathematical research field (and mainly concerning applications in Economics) and by M. Huang, P.E. Caines and R.P. Malham\'{e} \cite{HCM} for the engineering field.
It generalizes the mean field theory from Statistical Physics. Both theories deal with large population of individuals/ agents/ particles by letting their number tend to infinity. The main difference is that MFG theory considers each individual/ agent as a player capable of taking decisions in order to minimize or maximize 
an optimization criterion. The MFG theory studies large populations of interacting individuals/ agents, each one minimizing or maximizing an optimization criterion while taking into account the other players' actions. Each individual/ agent makes its decision (strategy) based on its access to the 
information about the structure of the population and based on the statistical information about the surrounding population. Actually, MFG theory shows that the strategies to be implemented are based on the distribution of the other players.
MFG pays a special attention to the equilibria in these games. 

The MFG theory has several applications in Population Dynamics, Economics, Social Networks, etc.
In Population Dynamics for example the MFG theory concerns the displacement of each individual of a large population in order to minimize an individual cost and taking into account the displacement of the other individuals. Each individual moves/acts (in a competitive environment) based on the probabilistic distribution 
of the individuals (probability density of the population) instead of considering the trajectories of one or of another individual. Each individual moves in order to minimize a mean cost (which is the Mean Field) which encodes the probabilistic distribution of the population.
It is of great interest in finding equilibria and the related Nash-type strategies.    

A standard MFG system consists of a Hamilton-Jacobi equation for a value function $p(t,x)$ and a Focker-Planck (FP) equation for the probabilistic distribution $\rho (t,x)$ of the individuals/ agents. 
\vspace{3mm}

Let us say a few words about the MFG system related to the following stochastic optimal control problem:
$$\underset{u \in {\cal U}}{\text{Minimize}} \ \mathbb{E}\Big[ \int_0^T[L(t,X(t),u(t,X(t)))+F(t,X(t), \rho (t,X(t)))] \, dt\Big]  \leqno({\bf P_S})$$
$$+ \mathbb{E}\Big[ F^0(X(T),\rho (T, X(T)))\Big], $$
where $T$ is a positive constant, $(X(t))_{t\in [0,T]}$ is the solution (or weak/ martingale solution) to the following stochastic differential equation
\begin{equation}\label{eqSDE}
\left\{ \begin{array}{ll}
\displaystyle dX(t)=f(t,X(t),u(t,X(t))) \, dt+\sigma (t, X(t)) \, dW(t), & t\in [0,T], \vspace{2mm} \\
X(0)=X_0 \ , & ~
                                     \end{array}
                          \right. 
\end{equation}
and $\rho(t)$ is the probability density of $X(t)$, for any $t\in[0,T]$.
Here $(\Omega , {\cal F}, \mathbb{P})$ is a probability space, $W:[0,T]\times \Omega \longrightarrow \mathbb{R}^N$ ($N\in \mathbb{N}^*$) is a Wiener process, and $({\cal F}_t)_{t\in [0,T]}$ is the corresponding natural filtration. Moreover, we
assume that $X_0$ is an $\mathbb{R}^d$-valued random variable ($d\in \mathbb{N}^*$) independent of the Wiener process and admits a density $\rho _0$. We have that
$f : [0,T]\times \mathbb{R}^d \times \mathbb{R}^{\ell}\longrightarrow {\cal M}_{d\times  \ell}(\mathbb{R})$ (${\ell }\in \mathbb{N}^*$), $\sigma : [0,T]\times \mathbb{R}^d\longrightarrow {\cal M}_{d\times  N}(\mathbb{R})$ are Borel functions. 
The set of controllers is
$${\cal U}= \left\{v : [0, T] \times \mathbb{R}^d\longrightarrow \mathbb{R}^{\ell }; \, v \ \mbox{\rm is a Borel function}, \ v(t,x)\in U_0 \ \mbox{\rm a.e. } (t,x)\in Q_T\right\} \, ,$$
where $U_0$ is a closed, bounded  and convex subset of $\mathbb{R}^{\ell }$ and $Q_T=(0,T)\times \mathbb{R}^d$. 
Let $a:[0,T]\times \mathbb{R}^d\longrightarrow {\cal M}_d(\mathbb{R})$ be defined by $a(t,x)=\frac{1}{2}\sigma (t,x)\sigma ^T(t,x)$. 
We have that $L:[0,T]\times \mathbb{R}^d\times \mathbb{R}^{\ell}\longrightarrow \mathbb{R}\cup \{+\infty\}$, $F:[0,T]\times \mathbb{R}^d\times \mathbb{R}\longrightarrow \mathbb{R}$, $F^0:\mathbb{R}^d\times \mathbb{R}\longrightarrow \mathbb{R}$.
\vspace{3mm}

$(P_S)$ models the problem of finding an optimal displacement  $(u^*(t,X(t)))_{t\in [0,T]}$ of an individual of a large population (with respect to the cost functional in $(P_S)$). 
Let us assume that $\rho(t)$ is a fixed probability density in $\mathbb{R}^d$, for any $t\in[0,T]$. So, problem $(P_S)$ may be viewed now as 
an optimal control problem depending on $\rho $.

Let us denote by
$$\psi ^{\rho }(t,x)= \underset{u \in {\cal U}_t}{\text{inf}} \ \Big\{ \mathbb{E}\Big[ \int_t^T[L(t,X(s),u(s,X(s)))+F(s,X(s), \rho (s,X(s)))] \, dt\Big] $$
$$+ \mathbb{E}\Big[ F^0(X(T),\rho (T, X(T)))\Big] \Big\} \, ,$$
where $(X(s))_{s\in [t,T]}$ is the solution to
\begin{equation*}
\left\{ \begin{array}{ll}
\displaystyle dX(s)=f(s,X(s),u(s,X(s))) \, dt+\sigma (s, X(s)) \, dW(s), & s\in [t,T], \vspace{2mm} \\
X(t)=x\in \mathbb{R}^d \, ,  & ~
                                     \end{array}
                          \right. 
\end{equation*}
and 
$${\cal U}_t=\left\{v : [t, T] \times \mathbb{R}^d\longrightarrow \mathbb{R}^{\ell }; \, v \ \mbox{\rm is a Borel function}, \ v(s,x)\in U_0 \ \mbox{\rm a.e. } (s,x)\in (t,T)\times \mathbb{R}^d \right\}.$$
Under appropriate hypotheses it follows that $\psi ^{\rho }$ is a generalized solution to the following Hamilton-Jacobi equation (see e.g. \cite{FR}):
\begin{equation*}
\left\{ \begin{array}{ll}
\displaystyle
\frac{\partial \psi }{\partial t}(t,x)-\max_{v\in U_0}\{ -f(t,x,v)\cdot \nabla \psi (t,x)-L(t,x,v)\} & \vspace{2mm} \\
\hspace{45mm} \displaystyle +a_{ij}(t,x)\frac{\partial ^2\psi }{\partial x_i\partial x_j}(t,x)+F(t,x,\rho (t,x))=0, \quad &(t,x)\in Q_T,\vspace{2mm} \\
\psi(T,x)=F^0(x,\rho (T,x)), & x\in \mathbb{R}^d \, 
                                     \end{array}
                          \right. 
\end{equation*}
and that $u^{\rho }(t,x)= \mbox{\rm argmax} \, \{ -f(t,x,v)\cdot \nabla \psi ^{\rho }(t,x)-L(t,x,v); \ v\in U_0\} $ is an optimal controller for $(P_S)$-(\ref{eqSDE}).

If we assume that $L$ is convex with respect to $u$ ($L(t,x,v)<+\infty $ for $v\in U_0$ and $L(t,x,v)=+\infty $ for $v\in \mathbb{R}^{\ell }\setminus U_0$), then
$$-f_u(t,x,u^{\rho }(t,x))^T\nabla \psi ^{\rho }(t,x)\in L_u(t,x,u^{\rho }(t,x)), \quad \mbox{\rm for } (t,x)\in Q_T \, ,$$
where $f_u(t,x,u)$ denotes the differential of $u\mapsto f(t,x,u)$ at $u$, $f_u(t,x,u)^T$ is the transpose of the matrix $f_u(t,x,u)$ and $L_u(t,x,u)$ is the subdifferential of $u\mapsto L(t,x,u)$ at $u$.
It follows that
$$-f(t,x,u^{\rho }(t,x))\cdot \nabla \psi ^{\rho }(t,x)-f_u(t,x,u^{\rho }(t,x))(v-u^{\rho }(t,x))\cdot \nabla \psi ^{\rho }(t,x)-L(t,x,v)$$
$$\leq -f(t,x,u^{\rho }(t,x))\cdot \nabla \psi ^{\rho }(t,x)-L(t,x,u^{\rho }(t,x)), \quad \mbox{\rm for } (t,x)\in Q_T, v\in \mathbb{R}^{\ell }$$
and consequently 
$$\max_{v\in U_0}\{ -f(t,x,v)\cdot \nabla \psi ^{\rho }(t,x)-L(t,x,v)\} $$
$$=\max_{v\in U_0}\{  -f(t,x,u^{\rho }(t,x))\cdot \nabla \psi ^{\rho }(t,x)-f_u(t,x,u^{\rho }(t,x))(v-u^{\rho }(t,x))\cdot \nabla \psi ^{\rho }(t,x)-L(t,x,v) \}$$
$$= -f(t,x,u^{\rho }(t,x))\cdot \nabla \psi ^{\rho }(t,x)+f_u(t,x,u^{\rho }(t,x))u^{\rho }(t,x)\cdot \nabla \psi ^{\rho }(t,x)$$
$$+\max_{v\in U_0}\{  -f_u(t,x,u^{\rho }(t,x))v\cdot \nabla \psi ^{\rho }(t,x)-L(t,x,v) \}$$
$$= -f(t,x,u^{\rho }(t,x))\cdot \nabla \psi ^{\rho }(t,x)+f_u(t,x,u^{\rho }(t,x))u^{\rho }(t,x)\cdot \nabla \psi ^{\rho }(t,x)$$
$$+H(t,x,  -f_u(t,x,u^{\rho }(t,x))^T\nabla \psi ^{\rho }(t,x)),$$
for  $(t,x)\in Q_T$. Here $H(t,x,\cdot )$ is the Hamiltonian corresponding to the Lagrangean $L(t,x,\cdot )$ for any $(t,x)\in Q_T$, i.e.
$$H(t,x,r)=\sup \ \left\{r\cdot v-L(t,x,v); \ v\in \mathbb{R}^{\ell } \right\} , \quad (t,x)\in Q_T $$ 
(see \cite{FR}). The Hamilton-Jacobi equation may be rewritten as
\begin{equation*}
\left\{ \begin{array}{ll}
\displaystyle
\frac{\partial \psi }{\partial t}(t,x)-H(t,x,-f_u(t,x,u^{\rho }(t,x))^T\nabla \psi (t,x)) & \vspace{2mm} \\
\hspace{15mm} \displaystyle +a_{ij}(t,x)\frac{\partial ^2\psi }{\partial x_i\partial x_j}(t,x)+F(t,x,\rho (t,x)) & \vspace{2mm} \\
\hspace{15mm} =-(f(t,x,u^{\rho }(t,x))-f_u(t,x,u^{\rho }(t,x))u^{\rho }(t,x))\cdot \nabla \psi (t,x), \ &(t,x)\in Q_T,\vspace{2mm} \\
\psi(T,x)=F^0(x,\rho (T,x)), & x\in \mathbb{R}^d \, .
                                     \end{array}
                          \right. 
\end{equation*}
On the other hand we formally get that
\begin{equation}\label{eq_uopt}
u^{\rho }(t,x)=H_q(t,x,-f_u(t,x,u^{\rho }(t,x))^T\nabla \psi ^{\rho }(t,x)), \quad (t,x)\in Q_T \, ,
\end{equation}
where $H_q(t,x,q)$ is the subdifferential of $q\mapsto H(t,x,q)$ at $q$.
If $\tilde{X}^{\rho }$ is the corresponding optimal process for $(P_S)$-(\ref{eqSDE}), with $u$ given by (\ref{eq_uopt}), then its probability density $\tilde{\rho }^{\rho }$ is a distributional solution and $t$-narrowly continuous in $L^1(\mathbb{R}^d)$ to the following Fokker-Planck equation
\begin{equation}\label{eqMFG'}
\left\{ \begin{array}{ll}
\displaystyle
\frac{\partial \tilde{\rho}}{\partial t}(t,x)+\nabla \cdot \left(f(t,x,u^{\rho }(t,x)) \tilde{\rho }(t,x)\right) -\frac{\partial ^2(a_{ij}\tilde{\rho })}{\partial x_i\partial x_j}(t,x)=0, \quad &(t,x)\in Q_T,\vspace{2mm} \\
\tilde{\rho }(0,x)=\rho_0(x), &x\in \mathbb{R}^d \, .
\end{array}
\right.
\end{equation}

Conversely, if $\tilde{\rho }^{\rho }$ is a distributional solution to (\ref{eqMFG'}), $t$-narrowly continuous in $L^1(\mathbb{R}^d)$ and $\tilde{\rho }^{\rho }(t)$ is a probability density for any $t\in[0,T]$, where $u(t,x):=u^{\rho}(t,x)$, $(t,x)\in Q_T$, then via superposition's principle (see \cite{figalli}, \cite{trevisan}, \cite{HRW}) we get that there exists $\big(\tilde{X}^{\rho }(t)\big)_{t\in[0,T]}$ a weak (martingale) solution to (\ref{eqSDE}) which admits $\tilde{\rho }^{\rho }(t)$ as a probability density for any $t\in[0,T]$.

To summarize: for a fixed $\rho $ we get an optimal pair $(\tilde{u},\tilde{\rho })$ for $(P_S)$-(\ref{eqSDE}). The mean field game strategy asks to find an anticipated $\rho $ such that
$\tilde{\rho }^{\rho }=\rho $. If $\tilde{\rho }=\rho $ and taking $p=-\psi ^{\rho }$ we get that $(\rho ,p)$ is a solution to 
\begin{equation}\label{eqMFG}
\left\{ \begin{array}{ll}
\displaystyle
\frac{\partial \rho}{\partial t}(t,x)+\nabla \cdot \left(f^u(t,x)\rho(t,x)\right) -\frac{\partial ^2(a_{ij}\rho )}{\partial x_i\partial x_j}(t,x)=0, \ &(t,x)\in Q_T,\vspace{2mm} \\
\displaystyle
\frac{\partial p}{\partial t}(t,x)+H\left(t,x,f_u^u(t,x)^T\nabla p(t,x)\right) +a_{ij}(t,x)\frac{\partial ^2p}{\partial x_i\partial x_j}(t,x) & \vspace{2mm} \\
\hspace{25mm} \displaystyle +(f^u(t,x)-f_u^u(t,x)u(t,x))\cdot \nabla p(t,x)=F(t,x,\rho (t,x)), \ &(t,x)\in Q_T,\vspace{2mm} \\
\rho(0,x)=\rho_0(x), \quad p(T,x)+F^0(x,\rho (T,x))=0, & x\in \mathbb{R}^d \, ,
\end{array}
\right.
\end{equation}
where $f^u(t,x)=f(t,x,u(t,x))$ for $(t,x)\in Q_T$ and $u(t,x)=H_q(t,x,f_u^u(t,x)^T\nabla p(t,x))$.
The corresponding strategy is a Nash-type equilibrium: Each individual chooses his optimal strategy taking into account the global information that is available to him and that
results from the actions of all individuals.
\vspace{3mm}

The present paper concerns some basic properties of the MFG system (\ref{eqMFG}), where $H(t,x,q)$ is  is a measurable function with respect to $(t,x)\in Q_T$ for any $q\in \mathbb{R}^{\ell}$ and convex and continuous with respect to $q\in \mathbb{R}^{\ell }$ for almost any $(t,x)\in Q_T$, $F(t,x,r)=G_r(t,x,r)$ and $F^0=G^0_r(x,r)$ for $(t,x,r)\in Q_T\times \mathbb{R}$,
and $(x,r)\in \mathbb{R}^{d+1}$, respectively, where $G$ and $G^0$ are measurable with respect to $(t,x)\in Q_T$ and $x\in \mathbb{R}^d$, respectively, and convex and continuous with respect to $r\in \mathbb{R}$. 
\vspace{3mm}

Let us recall that the existence and uniqueness of classical or weak solutions to (\ref{eqMFG}) has been investigated by several researchers; see e.g. \cite{LL}, \cite{porretta2015}, \cite{porretta2016}, \cite{CD}, \cite{cardaliaguet2017}, \cite{cardaliaguet2021}.
Their results mainly refer to the case when the domain is a $d$-dimensional torus and functions $F, F^0, H$ are smooth. The approaches use fixed point arguments and refined results for linear parabolic equations.
For other contributions we mention \cite{CGPT}, \cite{MS2015}, \cite{MS2016}, \cite{OS}. We notice that we use in our present paper a different approach and that our hypotheses on $f, \sigma , F, F^0, H$ are less restrictive.

For MFG in infinite dimension we refer to \cite{SGS2024}.
\vspace{3mm}

The approach in the present paper employs an idea previously used in \cite{LL} and \cite{barbu2024} (for different situations) which is based on the remark that (\ref{eqMFG}) is the optimality system corresponding to the next deterministic optimal control problem.

$$\underset{u \in {\cal U}}{\text{Minimize}} \ \int_{Q_T}[L(t,x,u(t,x))\rho (t,x)+G(t,x, \rho (t,x))] \, dx \, dt  + \int_{\mathbb{R}^d} G^0(x,\rho (T, x)) \, dx \, , \leqno({\bf P})$$
subject to $\rho $, weak solution to the following Fokker-Planck equation
\begin{equation}\label{eqFP}
\left\{ \begin{array}{ll}
\displaystyle
\frac{\partial \rho}{\partial t}(t,x)=-\nabla \cdot (f^u(t,x)\rho(t,x))+\frac{\partial ^2(a_{ij}\rho )}{\partial x_i\partial x_j}(t,x), \quad &(t,x)\in Q_T,\vspace{2mm} \\
\rho(0,x)=\rho_0(x), & x\in \mathbb{R}^d \, ,
\end{array}
\right.
\end{equation}
where $G:Q_T\times \mathbb{R}\longrightarrow \mathbb{R}$, $G^0:\mathbb{R}^{d+1}\longrightarrow \mathbb{R}$ are such that
$$F(t,x,r)=G_r(t,x,r), \quad F^0(x,r)=G^0_r(x,r), \quad \forall t\in (0,T), \, x\in \mathbb{R}^d, \, r\in \mathbb{R} $$
and $L$ is given by 
$$L(t,x,u)=\sup \ \left\{ u\cdot r-H(t,x,r); \, r\in \mathbb{R}^{\ell }\right\}, \quad (t,x)\in Q_T \, .$$ 
\vspace{2mm}

Here is the structure of the paper. In Section 2 we establish some basic properties of the weak solution to a linear Fokker-Planck equation (where $f^u(t,x)$ is replaced by a general $\tilde{f}(t,x)$).
The existence of a weak solution for the MFG system is established in Section 3. A uniqueness for the weak solution of the MFG system is treated in Section 4. Some comments and further investigation directions are 
given in Section 5.

\vspace{3mm}

\noindent
\bf Conventions and notations. \rm We use the Einstein summation convention and the following notations
\begin{itemize}
\item{} ``$\cdot $'' for the usual scalar product of $\mathbb{R}^n$, $|\cdot |$ for the Euclidean norm of $\mathbb{R}^n$;
\item{} $L^k$ for $L^k(Q)$ or $L^k(Q;\mathbb{R}^n)$ and $\| \cdot \|_{L^k}$, when there is no danger of confusion and $H^1$ for $H^1(\mathbb{R}^d)$ and $H^{-1}$ for $H^{-1}(\mathbb{R}^d)$;
\item{} $\langle \cdot ,\cdot \rangle _{L^2}$ for the usual scalar product in $L^2$ and $\langle \cdot ,\rangle _{H^{-1},H^1}$ for the duality between $H^{-1}$ and $H^1$;
\item{} $f^u(t,x)$ for $f(t,x,u)$;
\item{} $A^T$ for the transpose of the matrix $A$;
\item{} $f_u(t,x,u)$ or $f_u^u(t,x)$ for the differential with respect to $u$ of $u\mapsto f(t,x,u)$ and $u\mapsto f^u(t,x)$, respectively;
\item{} $G_r(t,x,r)$ and $G^0_r(x,r)$ for the differential or subdifferential with respect to $r$ of of $r\mapsto G(t,x,r)$ and $r\mapsto G^0(x,r)$, respectively;
\item{} $L_u(t,x,u)$ for the subdifferential with respect to $u$ of $u\mapsto L(t,x,u)$,
\item{} $H_q(t,x,q)$ for the differential or subdifferential (which may be multivalued) with respect to $q$ of $q\mapsto H(t,x,q)$;
\item{} $\rho^u$ for the weak solution to (\ref{eqFP}).
\end{itemize}
We shall denote by $c, C, \tilde{c}, \tilde{C}, c_k, C_k, ...$, several positive constants occurring in the present paper.

\section{Analysis of a linear Fokker-Planck equation. The case when $\rho _0\in L^1(\mathbb{R}^d)\cap L^2(\mathbb{R}^d)$}

Consider the following Fokker-Planck equation

\begin{equation}\label{eqLFP}
\left\{ \begin{array}{ll}
\displaystyle
\frac{\partial \rho}{\partial t}(t,x)=-\nabla \cdot (\tilde{f}(t,x)\rho(t,x))+\frac{\partial ^2(a_{ij}\rho )}{\partial x_i\partial x_j}(t,x), \quad &(t,x)\in Q_T,\vspace{2mm} \\
\rho(0,x)=\rho_0(x), & x\in \mathbb{R}^d \, ,
\end{array}
\right.
\end{equation}

Assume that the following hypotheses hold:
\begin{itemize}
\item[\bf (H0)] $\tilde{f}\in L^{\infty }(Q_T;\mathbb{R}^d)$;
\item[\bf (H1)] $\rho _0\in L^1(\mathbb{R}^d)\cap L^2(\mathbb{R}^d)$;
\item[\bf (H2)] $\sigma \in C^{1,3}_b([0,T]\times \mathbb{R}^d; \mathbb{R}^{d\cdot N})$ and there exists a positive constant $\gamma $ such that
$$\sigma (t,x)\sigma (t,x)^Ty\cdot y=a_{ij}(t,x)y_iy_j\geq \gamma |y|^2, \quad \forall (t,x,y)\in [0,T]\times \mathbb{R}^d\times \mathbb{R}^d \, .$$
Here $a_{ij}(t,x)=\sigma_{ik}(t,x)\sigma_{jk}(t,x), \quad \forall (t,x)\in [0,T]\times \mathbb{R}^d, \ i,j \in \{ 1,2,...,d\} $.
\end{itemize}

\noindent
\bf Definition 2.1. \rm We say that $\rho $ is a \it weak solution \rm to (\ref{eqLFP}) if
\begin{equation}\label{eq3}
\left\{ \begin{array}{ll}
\rho \in L^2(0,T; H^1(\mathbb{R}^d))\cap W^{1,2}([0,T];H^{-1}(\mathbb{R}^d)), & \vspace{2mm} \\
\displaystyle
\frac{d \rho}{d t}(t)=-\nabla \cdot (\tilde{f}(t)\rho(t))+\frac{\partial ^2(a_{ij}(t)\rho (t)) }{\partial x_i\partial x_j} \quad \mbox{\rm in }H^{-1}(\mathbb{R}^d), & \mbox{\rm a.e. } t\in (0,T) ,\vspace{2mm} \\
\rho(0,x)=\rho_0(x), & \mbox{\rm a.e. } x\in \mathbb{R}^d \, .
\end{array}
\right.
\end{equation}
\vspace{2mm}

Let us recall that if $\rho \in L^2(0,T; H^1(\mathbb{R}^d))\cap W^{1,2}([0,T];H^{-1}(\mathbb{R}^d))$, then $\rho \in C([0,T]; L^2(\mathbb{R}^d))$ and that
\begin{equation*}
\frac{d}{dt}\langle \rho (t),\varphi\rangle _{L^2}=\Big\langle \frac{d\rho }{dt}(t),\varphi \Big\rangle _{H^{-1},H^1}, \quad \mbox{\rm a.e. } t\in (0,T), \ \forall \varphi \in H^1(\mathbb{R}^d) \, .
\end{equation*}

\noindent
\bf Theorem 2.1. \it Problem (\ref{eqLFP}) admits a unique weak solution $\rho $. Moreover, there exists a positive constant $C_T$ (depending on $T$ and $\|\tilde{f}\|_{L^{\infty }}$) such that
\begin{equation}\label{eq5}
\|\rho (t)\|_{L^2}^2+\int_0^t\|\nabla \rho (s)\|_{L^2}^2ds\leq C_T^2\|\rho _0\|^2_{L^2}, \quad \forall t\in [0,T] \, ,
\end{equation}
\begin{equation}\label{eq6}
\|\rho (t)\|_{L^m}\leq C_T\|\rho _0\|_{L^m}, \quad \forall t\in [0,T], \ \forall m\in [1, +\infty ] \, .
\end{equation}

If $\rho _0(x)\geq 0$ a.e. $x\in \mathbb{R}^d$, then $\rho (t,x)\geq 0$, $\forall t\in [0,T]$, a.e. $x\in \mathbb{R}^d$. If 
$$\rho _0\in {\cal P}=\left\{ \psi \in L^1(\mathbb{R}^d); \ \psi (x)\geq 0 \ a.e. \ x\in \mathbb{R}^d, \ \int_{\mathbb{R}^d}\psi (x)dx=1 \right\} ,$$
then $\rho (t)\in {\cal P}$, $\forall t\in [0,T]$.

If $\rho _0(x)> 0$ a.e. $x\in \mathbb{R}^d$ and $\log \rho _0\in L^1_{loc}(\mathbb{R}^d)$, then 
\begin{equation*}
\begin{array}{ll}
\rho (t,x)>0, \quad a.e. \ (t,x)\in Q_T, \vspace{2mm} \\
\log \rho (t)\in L^1_{loc}(\mathbb{R}^d), \quad \forall t\in [0,T] \, .
\end{array}
\end{equation*}

Proof. \rm
(\ref{eqLFP}) may be written in the abstract form
\begin{equation*}
\left\{ \begin{array}{ll}
\displaystyle
\frac{d\rho }{dt}(t)+A(t)\rho (t)=0, \quad t\in (0,T) \, ,\vspace{2mm} \\
\rho (0)=\rho _0 \, ,
\end{array}
\right.
\end{equation*}
where 
$$\Big\langle A(t)\varphi , \psi \Big\rangle_{H^{-1},H^1}=-\int_{\mathbb{R}^d}\tilde{f}(t,x)\varphi (x)\cdot \nabla \psi (x) \, dx+\int_{\mathbb{R}^d}\frac{\partial }{\partial x_j}(a_{ij}(t,x)\varphi (x))  \frac{\partial \psi }{\partial x_i}(x) \, dx $$
$$={\cal A}(t;\varphi ,\psi ) \, ,$$
$\forall \varphi, \psi \in H ^1(\mathbb{R} ^d)$.
We have that 
\begin{itemize}
\item[(j)] ${\cal A}(t; \cdot ,\cdot )$ is bilinear a.e. $t\in (0,T)$, 
$$t\mapsto {\cal A}(t;\varphi ,\psi ) \quad \mbox{\rm is measurable for any } \varphi ,\psi \in H^1(\mathbb{R}^d) \, ;$$
\item[(jj)] There exists a positive constant $M$ such that
$$|{\cal A}(t; \varphi ,\psi )|\leq M\| \varphi \|_{H^1}\| \psi \| _{H^1}, \quad \mbox{\rm a.e. } t\in (0,T), \ \forall \varphi ,\psi \in H^1(\mathbb{R}^d) \, ;$$
\item[(jjj)] There exist $\alpha ^*, \beta ^*\in \mathbb{R}$, $\alpha ^*>0$ such that
$${\cal A}(t; \varphi ,\varphi )\geq \alpha ^*\| \varphi \|^2_{H^1}-\beta ^*\| \varphi \|^2_{L^2}, \quad \mbox{\rm a.e. } t\in (0,T), \ \forall \varphi \in H^1 \, .$$
\end{itemize}
By Lions' existence theorem (see \cite{brezis}) we conclude that there exists a unique $\rho \in L^2(0,T; H^1(\mathbb{R}^d))\cap W^{1,2}([0,T];H^{-1}(\mathbb{R}^d))$ such that
$$\Big\langle \frac{d\rho }{dt}(t), \varphi \Big\rangle _{H^{-1},H^1}+{\cal A}(t;\rho (t),\varphi )=0, \quad \mbox{\rm a.e. } t\in (0,T), \ \forall \varphi \in H^1(\mathbb{R}^d) $$
and $\rho (0)=\rho _0$.
It means that (\ref{eqLFP}) admits a unique weak solution.
\vspace{3mm}

If we assume now that $\rho _0(x)\geq 0$ a.e., then we multiply (\ref{eqLFP}) by $\rho ^-(t)$ and we get that
$$\rho (t,x)\geq 0, \quad \mbox{\rm a.e. } x\in \mathbb{R}^d, \ \forall t\in [0,T] \, . $$
If moreover, $\rho _0\in {\cal P}$, then it follows (as in \cite{anita2021}) that $\rho (t)\in {\cal P}$, for any $t\in [0,T]$.
\vspace{3mm}

Assume now that $\rho _0(x)>0$ a.e. $x\in \mathbb{R}^d$ and $\log \rho _0\in L^1_{loc}(\mathbb{R}^d)$. Let 
$$\psi _0(x)=\left\{ \begin{array}{ll}
\displaystyle
e^{\frac{1}{|x|^2-1}}, \quad & x\in \mathbb{R}^d, \ |x|<1, \vspace{2mm} \\
0, &  x\in \mathbb{R}^d, \ |x|\geq 1 \, .
\end{array}
\right. $$  
It is well known that $\psi _0\in C_0^{\infty }(\mathbb{R}^d)$, $supp \ \psi _0=\overline{B(0_d;1)}$.

For an arbitrary but fixed $R>0$ we have that $\psi (x)=\psi _0(\frac{x}{R})$ satisfies
$\psi \in C_0^{\infty }(\mathbb{R}^d)$ and $supp \ \psi =\overline{B(0_d;R)}$. On the other hand there exists a positive constant $M_R$ such that
$$| \nabla \psi (x)|^2\leq M_R\psi (x), \quad \forall x\in \mathbb{R}^d \, .$$

Let $\varepsilon >0$ be arbitrary but fixed and let $\varphi (t,x)=\psi (x)(\rho (t,x)+\varepsilon )^{-1}$.
If we multiply (\ref{eqLFP}) by $\varphi (t)$ we get that
$$\frac{d}{dt}\int_{\mathbb{R}^d}\psi (x)\log (\rho (t,x)+\varepsilon )dx +\int_{\mathbb{R}^d}\frac{\partial }{\partial x_j}(a_{ij}(t,x)\rho (t,x))\frac{\partial }{\partial x_i}\Big( \frac{\psi (x)}{\rho (t,x)+\varepsilon}\Big)dx $$
$$=\int_{\mathbb{R}^d}{\tilde f}(t,x)\rho (t,x)\cdot \nabla \Big( \frac{\psi (x)}{\rho (t,x)+\varepsilon }\Big) dx \, , \quad \mbox{\rm a.e. } t\in (0,T) \, .$$
This implies that
\begin{equation*}
\begin{array}{ll}
\displaystyle
\frac{d}{dt}\int_{\mathbb{R}^d}\psi (x)\log (\rho (t,x)+\varepsilon )dx \vspace{2mm} \\
\displaystyle
 \ \ \ \ \ +\int_{\mathbb{R}^d}\Big[ a_{ij}(t,x)\frac{\partial \rho }{\partial x_j}(t,x)+\frac{\partial a_{ij}}{\partial x_j}(t,x)\rho(t,x) \Big] \Big[ \frac{\frac{\partial \psi }{\partial x_i}(x)}{\rho (t,x)+\varepsilon }
-\frac{\psi (x)\frac{\partial \rho }{\partial x_i}(t,x)}{(\rho (t,x)+\varepsilon )^2}\Big] dx \vspace{2mm} \\
\displaystyle 
 \ \ \ =\int_{\mathbb{R}^d}{\tilde f}(t,x)\rho (t,x)\cdot \Big[ \frac{\nabla \psi (x)}{\rho (t,x)+\varepsilon }-\frac{\psi (x)\nabla \rho (t,x)}{(\rho (t,x)+\varepsilon )^2}\Big] dx \, , \quad \mbox{\rm a.e. } t\in (0,T) \, .
\end{array}
\end{equation*}
It follows that for almost any $t\in (0,T)$:
\begin{equation*}
\begin{array}{ll}
\displaystyle
\frac{d}{dt}\int_{\mathbb{R}^d}\psi \log (\rho (t)+\varepsilon )dx -\int_{\mathbb{R}^d}a_{ij}(t)\psi \frac{\frac{\partial \rho }{\partial x_i}(t)\frac{\partial \rho }{\partial x_j}(t)}{(\rho (t)+\varepsilon )^2}dx \vspace{2mm} \\
\displaystyle
=-\int_{\mathbb{R}^d}a_{ij}(t)\frac{\partial \rho }{\partial x_j}(t)\frac{\partial \psi }{\partial x_i}\frac{1}{\rho (t)+\varepsilon }dx
-\int_{\mathbb{R}^d}\frac{\partial a_{ij}}{\partial x_j}(t)\frac{\partial \psi }{\partial x_i}\frac{\rho (t)}{\rho (t)+\varepsilon }dx \vspace{2mm} \\
\displaystyle 
+\int_{\mathbb{R}^d}\frac{\partial a_{ij}}{\partial x_j}(t)\psi \rho (t)\frac{\frac{\partial \rho }{\partial x_i}(t)}{(\rho (t)+\varepsilon )^2}dx 
+\int_{\mathbb{R}^d}{\tilde f}(t)\rho (t)\cdot \Big[ \frac{\nabla \psi }{\rho (t)+\varepsilon }-\frac{\psi \nabla \rho (t)}{(\rho (t)+\varepsilon )^2}\Big] dx 
\end{array}
\end{equation*}
(for the sake of simplicity we have omitted everywhere the argument $x$).
If we integrate on $[0,t]$ we get that there exist $c_1, c_2, c_3, c_4, c_5$, positive constants such that for any $t\in [0,T]$:
\begin{equation*}
\begin{array}{ll}
\displaystyle
-\int_{\mathbb{R}^d}\psi (x)\log (\rho (t,x)+\varepsilon )dx +\int_{\mathbb{R}^d}\psi (x)\log (\rho _0(x)+\varepsilon)dx
+\gamma \int_0^t\int_{\mathbb{R}^d}\psi (x)\frac{|\nabla \rho (s,x)|^2}{(\rho (s,x)+\varepsilon )^2}dx \ ds \vspace{2mm} \\
\displaystyle 
\leq c_1\int_0^t\int_{\mathbb{R}^d}\frac{|\nabla \psi | |\nabla \rho |}{\rho +\varepsilon }dx \ ds 
+c_2\int_0^t\int_{\mathbb{R}^d}|\nabla \psi | \frac{\rho }{\rho +\varepsilon }dx \ ds
+c_3\int_0^t\int_{\mathbb{R}^d}\psi \frac{\rho }{\rho +\varepsilon }\frac{|\nabla \rho |}{\rho +\varepsilon }dx \ ds \vspace{2mm} \\
\displaystyle
 \ \ +c_4\int_0^t\int_{\mathbb{R}^d}|\nabla \psi | \frac{\rho }{\rho +\varepsilon }dx \ ds
+c_5\int_0^t\int_{\mathbb{R}^d}\psi \frac{\rho }{\rho +\varepsilon }\frac{|\nabla \rho |}{\rho +\varepsilon }dx \ ds \, .
\end{array}
\end{equation*}
After an easy calculation (and using the fact that $|\nabla \psi (x)| ^2\leq M_R\psi (x)$) we get that there exists a positive constant $c_6$ (depending on $R$) such that for any $t\in [0,T]$:
\begin{equation*}
\begin{array}{ll}
\displaystyle
-\int_{\mathbb{R}^d}\psi (x)\log (\rho (t,x)+\varepsilon )dx +\int_{\mathbb{R}^d}\psi (x)\log (\rho _0(x)+\varepsilon)dx
+\gamma \int_0^t\int_{\mathbb{R}^d}\psi (x)\frac{|\nabla \rho (s,x)| ^2}{(\rho (s,x)+\varepsilon )^2}dx \ ds \vspace{2mm} \\
\displaystyle 
\leq c_6 \, t  \Big[ 1 +\| \psi \| _{L^1}+\| \nabla \psi \| _{L^1}\Big] \, 
\end{array}
\end{equation*}
and consequently
\begin{equation*}
\begin{array}{ll}
\displaystyle
-\int_{\mathbb{R}^d}\psi (x)\log (\rho (t,x)+\varepsilon )dx \vspace{2mm} \\
\displaystyle 
\leq-\int_{\mathbb{R}^d}\psi (x)\log (\rho _0(x)+\varepsilon)dx
+c_7 \, \Big[ 1 +\| \psi \| _{L^1}+\| \nabla \psi \| _{L^1}\Big] \, ,
\end{array}
\end{equation*}
where $c_7$ is a positive constant. The last inequality implies that
\begin{equation*}
\begin{array}{ll}
\displaystyle
-\int_{\mathbb{R}^d, \rho (t,x)+\varepsilon <1}\psi (x)\log (\rho (t,x)+\varepsilon )dx\leq \int_{\mathbb{R}^d, \rho (t,x)+\varepsilon \geq 1}\psi (x)\log (\rho (t,x)+\varepsilon )dx \vspace{2mm} \\
\displaystyle 
-\int_{\mathbb{R}^d}\psi (x)\log (\rho _0(x)+\varepsilon)dx
+c_7 \, \Big[ 1 +\| \psi \| _{L^1}+\| \nabla \psi \| _{L^1}\Big] 
\end{array}
\end{equation*}
and consequently
\begin{equation*}
\begin{array}{ll}
\displaystyle
\int_{\mathbb{R}^d}\psi (x) \, |\log (\rho (t,x)+\varepsilon )|dx\leq 2 \int_{\mathbb{R}^d, \rho (t,x)+\varepsilon \geq 1}\psi (x)\log (\rho (t,x)+\varepsilon )dx \vspace{2mm} \\
\displaystyle 
-\int_{\mathbb{R}^d}\psi (x) \, \log (\rho _0(x)) dx +c_7 \, \Big[ 1 +\| \psi \| _{L^1}+\| \nabla \psi \| _{L^1}\Big] \, .
\end{array}
\end{equation*}
It follows that for any $\varepsilon \in (0,1)$:
\begin{equation*}
\begin{array}{ll}
\displaystyle
\int_{\mathbb{R}^d}\psi (x) \, |\log (\rho (t,x)+\varepsilon )|dx\leq 2 \int_{\mathbb{R}^d,\rho (t,x)+\varepsilon \geq 1}\psi (x)| \rho (t,x)+\varepsilon -1| \, dx \vspace{2mm} \\
\displaystyle 
-\int_{\mathbb{R}^d}\psi (x) \, \log (\rho _0(x)) dx +c_7 \, \Big[ 1 +\| \psi \| _{L^1}+\| \nabla \psi \| _{L^1}\Big] 
\end{array}
\end{equation*}
and since $\rho \in C([0,T]; L^2(\mathbb{R}^d))$, there exists a positive constant $c_8$ such that
$$\int_{\mathbb{R}^d}\psi (x) \, |\log (\rho (t,x)+\varepsilon )|dx\leq c_8$$
and so by Fatou's lemma we get that
$$\int_{\mathbb{R}^d}\psi (x) \, |\log \rho (t,x)|dx\leq c_8 $$
and consequently $\log \rho (t) \in L^1(B(0_d;\frac{R}{2}))$.
Since $R$ is an arbitrary positive constant we get that
$$\log \rho (t)\in L^1_{loc}(\mathbb{R}^d), \quad \forall t\in [0,T] \, .$$
\vspace{3mm}

If we multiply (\ref{eqLFP}) (or (\ref{eq3})) by $\rho (t)$ we obtain after an easy calculation inequality (\ref{eq5}).
\vspace{3mm}

Let us now prove (\ref{eq6}). We first prove this result in the particular case when 
${\tilde f}\in C^{1,2}_b([0,T]\times \mathbb{R}^d)$ and $\rho _0\in C_0(\mathbb{R}^d)$.

Let $\Gamma (t,x,s,\xi )$ be the fundamental solution of
$$v\mapsto \frac{\partial v}{\partial t}-\frac{\partial ^2(a_{ij}v)}{\partial x_i\partial x_j} \,$$
and let $\tilde{\Gamma }(t,x,s,\xi )$ be the fundamental solution of 
$$v\mapsto \frac{\partial v}{\partial t}+\nabla \cdot ({\tilde f}v)-\frac{\partial ^2(a_{ij}v)}{\partial x_i\partial x_j} \, . $$
\vspace{2mm}

\noindent
(\ref{eqLFP}) may be rewritten as
\begin{equation*}
\left\{ \begin{array}{ll}
\displaystyle
\frac{\partial \rho}{\partial t}(t,x)-\frac{\partial ^2(a_{ij}\rho )}{\partial x_i\partial x_j}(t,x)={\tilde F}(t,x)=-\nabla \cdot ({\tilde f}(t,x)\rho (t,x)), \quad & (t,x)\in Q_T, \vspace{2mm} \\
\rho(0,x)=\rho_0(x), & x\in \mathbb{R}^d 
\end{array}
\right.
\end{equation*}
and let us prove that its weak solution satisfies
\begin{equation}\label{eq9}
\begin{array}{ll}
\displaystyle
\rho (t,x)&
\displaystyle 
=\int_{\mathbb{R}^d}\Gamma (t,x,0,\xi )\rho _0(\xi )d\xi +\int_0^t\int_{\mathbb{R}^d}\Gamma (t,x,s,\xi ){\tilde F}(s,\xi )d\xi \ ds \vspace{2mm} \\
~ & \displaystyle 
=\int_{\mathbb{R}^d}\Gamma (t,x,0,\xi )\rho _0(\xi )d\xi +\int_0^t\int_{\mathbb{R}^d}\nabla _{\xi }\Gamma (t,x,s,\xi )\cdot {\tilde f}(s,\xi )\rho (s,\xi )d\xi \ ds, \mbox{\rm a.e. } (t,x)\in Q_T.
\end{array}
\end{equation} 

Indeed, let us recall (see \cite{friedman} Theorem 4.5, p. 141-142, Theorem 4.7, p. 143) that there exist two positive constants $c, C$ and two constants, depending on $\tilde {f}$, $\tilde{c}, \tilde{C}$,  such that
$$\Big| \frac{\partial ^r\Gamma }{\partial x_i^r}(t,x,s,\xi )\Big| , \ \Big| \frac{\partial ^r\Gamma }{\partial \xi _i^r}(t,x,s,\xi )\Big|  \leq C (t-s)^{-\frac{d+r}{2}}\exp \Big( -c\frac{|x-\xi |^2}{t-s}\Big) ,$$
and
$$\Big| \frac{\partial ^r\tilde{\Gamma} }{\partial x_i^r}(t,x,s,\xi )\Big| , \ \Big| \frac{\partial ^r\tilde{\Gamma }}{\partial \xi _i^r}(t,x,s,\xi )\Big|  \leq \bar{C} (t-s)^{-\frac{d+r}{2}}\exp \Big( -\bar{c}\frac{|x-\xi |^2}{t-s}\Big) ,$$
for any $0\leq s <t\leq T$, $x, \xi \in \mathbb{R}^d$, $r\in \{ 0,1\} $.
\vspace{2mm}

There exists a classical solution $\rho $ to (\ref{eqLFP}) ($\rho \in C^{1,2}((0,T]\times \mathbb{R}^d)$, satisfies the PDE in (\ref{eqLFP}) for $(t,x)\in (0,T]\times \mathbb{R}^d$ and $\lim_{t\rightarrow 0+}\rho (t,x)=\rho _0(x)$ for any $x\in \mathbb{R}^d$) given by
$$\rho (t,x)=\int_{\mathbb{R}^d}\tilde{\Gamma }(t,x,0,\xi )\rho _0(\xi )d\xi \, .$$
The properties of $\tilde{\Gamma }$ and $\rho _0$ imply that $\rho \in L^{\infty }(Q_T)$.

Let us prove now that $\rho $ is also a weak solution to (\ref{eqLFP}). Let $t\in [0,T]$ be arbitrary but fixed. By Young's inequality for convolution products we get that:
$$\| \rho (t) \|_{L^2} \leq \tilde{C} \| \rho_0 \| _{L^2} \int_{\mathbb{R}^d}\frac{1}{t^{\frac{d}{2}}}e^{-\tilde{c}\frac{|x|^2}{t}}dx= \tilde{C_1} \| \rho_0 \| _{L^2},$$
where $\tilde{C_1}=\tilde{C} \theta _0$, and $\theta _0 =\int_{\mathbb{R}^d}e^{-\tilde{c}|z|^2}dz$. It follows that $\rho \in L^{\infty }(0,T;L^2(\mathbb{R}^d))$. 

Let us prove now that $\rho \in L^2(0,T; H^1(\mathbb{R}^d))$. Let $\zeta \in C^1(\mathbb{R}^d)$ be such that 
$$\begin{array}{llll}
\zeta (x)=1, \quad &\mbox{\rm if } | x|\leq 1 , \quad &\zeta (x)=0, \quad &\mbox{\rm if } | x|\geq 2 \\
0<\zeta (x)< 1, \quad &\mbox{\rm if }  1<| x|<2, \quad &\displaystyle \frac{| \nabla \zeta (x)|}{\sqrt{\zeta (x)}}\leq \overline{M}_0, \quad &\mbox{\rm if } | x|<2 
   \end{array}$$
(here $\overline{M}_0\geq 0$ is a constant). An example of such a function $\zeta $ is given in \cite{anita2022}.

For any $R>0$ we consider the function $\zeta _R:\mathbb{R}^d\rightarrow \mathbb{R}$,
$\zeta _R(x)=\zeta \left( \frac{1}{R}x\right)$, $x\in \mathbb{R}^d.$ It is obvious that $\zeta _R\in C^1(\mathbb{R}^d)$ and that
$$\begin{array}{llll}
\zeta _R(x)=1, \quad &\mbox{\rm if } | x|\leq R, \quad &\zeta _R(x)=0, \quad &\mbox{\rm if } | x|\geq 2R \\
0<\zeta _R(x)< 1, \quad &\mbox{\rm if } R<| x|<2R, \quad &\displaystyle \frac{| \nabla \zeta _R(x)|}{\sqrt{\zeta _R(x)}}\leq \frac{\overline{M}_0}{R}, \quad &\mbox{\rm if } | x|<2R . 
   \end{array}$$

By multiplying the Fokker-Planck equation in (\ref{eqLFP}) with $\rho \zeta_R$ (for $R\geq 1$) and integrating over $[0,t]\times \mathbb{R}^d$ we get the following:
\begin{equation}\label{eqa.1}
\begin{array}{cc}
\displaystyle
\int_0^t\int_{\mathbb{R}^d}\frac{\partial \rho }{\partial t}(s,x)\rho (s,x)\zeta_R(x)dx \ ds=-\int_0^t\int_{\mathbb{R}^d}\nabla \cdot ({\tilde f}(s,x)\rho (s,x))\rho (s,x)\zeta_R(x)dx \ ds \\
\displaystyle
+\int_0^t\int_{\mathbb{R}^d}\frac{\partial^2(a_{ij}\rho )}{\partial x_i \partial x_j}(s,x)\rho (s,x)\zeta_R(x)dx \ ds. 
\end{array}
\end{equation}

Since, $\frac{\partial \rho }{\partial t} \, \rho \, \zeta_R=\frac{1}{2}\frac{\partial }{\partial t}(\rho ^2 \zeta_R)$, it follows that for any $t\in [0,T]$:
$$\int_0^t\int_{\mathbb{R}^d}\frac{\partial \rho }{\partial t}(s,x)\rho (s,x)\zeta_R(x)dx \ ds=\frac{1}{2}\int_{\mathbb{R}^d}\rho (t,x)^2\zeta_R(x)dx-\frac{1}{2}\int_{\mathbb{R}^d}\rho (0,x)^2\zeta_R(x)dx$$
$$=\frac{1}{2}\int_{\mathbb{R}^d}\rho (t,x)^2\zeta_R(x)dx-\frac{1}{2}\int_{\mathbb{R}^d}\rho_0(x)^2\zeta_R(x)dx.$$
By using the divergence theorem we rewrite the right hand side term of (\ref{eqa.1}) as:
$$-\int_0^t\int_{\mathbb{R}^d}\nabla \cdot ({\tilde f}\rho )\rho \zeta_R \, dx \ ds+\int_0^t\int_{\mathbb{R}^d}\frac{\partial^2(a_{ij}\rho )}{\partial x_i \partial x_j}\rho \zeta_R \, dx \ ds$$
$$=\int_0^t\int_{\mathbb{R}^d}\rho {\tilde f}\cdot \nabla (\rho \zeta_R) \, dx \ ds-\int_0^t\int_{\mathbb{R}^d}\frac{\partial (a_{ij}\rho )}{\partial x_j}\frac{\partial (\rho \zeta_R)}{\partial x_i} \, dx \ ds$$
and by simple calculation
$$=\int_0^t\int_{\mathbb{R}^d}\rho {\tilde f}\cdot \nabla \rho \zeta_R \, dx \ ds+\int_0^t\int_{\mathbb{R}^d}\rho ^2{\tilde f}\cdot \nabla \zeta_R \, dx \ ds
-\int_0^t\int_{\mathbb{R}^d}\frac{\partial a_{ij}}{\partial x_j}\rho \frac{\partial \rho }{\partial x_i}\zeta_R \, dx \ ds$$
$$-\int_0^t\int_{\mathbb{R}^d}\frac{\partial a_{ij}}{\partial x_j}\rho ^2\frac{\partial \zeta_R}{\partial x_i} \, dx \ ds-\int_0^t\int_{\mathbb{R}^d}a_{ij}\frac{\partial \rho }{\partial x_j}\frac{\partial \rho }{\partial x_i}\zeta_R \, dx \ ds
-\int_0^t\int_{\mathbb{R}^d}a_{ij}\frac{\partial \rho }{\partial x_j}\rho \frac{\partial \zeta_R}{\partial x_i} \, dx \ ds.$$

After an easy calculation (and using the properties of ${\tilde f}$ and $\zeta _R$ and (H2)) we obtain that there exits $\overline{M}_1, \overline{M}_2, \overline{M}_3>0$ such that
$$\int_{\mathbb{R}^d}\rho (t,x)^2\zeta _R(x)dx +\gamma \int_0^t\int_{\mathbb{R}^d}|\nabla \rho (s,x)|^2\zeta _R(x)dx \ ds \leq \int _{\mathbb{R}^d}\rho _0(x)^2dx $$
$$+ \overline{M}_1\int_0^t\int_{\mathbb{R}^d}\rho (s,x)^2\zeta _R(x)dx \ ds +\overline{M}_2\int_0^t\int_{\mathbb{R}^d}\rho (s,x)^2dx \ ds \leq \overline{M}_3,$$
for any $t\in [0,T]$ and $R\geq 1$ (we have used the boundedness of $\rho $ and the properties of $\zeta _R$). It follows that
$$\int_{|x|\leq R}\rho (t,x)^2dx +\gamma \int_0^t\int_{|x|\leq R}|\nabla \rho (s,x)|^2dx \ ds \leq \overline{M}_3, \quad \forall t\in [0,T].$$
If we pass to the limit ($R\rightarrow +\infty $) we obtain that
$$\int_{\mathbb{R}^d}\rho (t,x)^2dx +\gamma \int_0^t\int_{\mathbb{R}^d}|\nabla \rho (s,x)|^2dx \ ds \leq \overline{M}_3, \quad \forall t\in [0,T],$$
and consequently that $\rho \in L^2(0,T;H^1(\mathbb{R}^d))$.
\vspace{3mm}

The next goal is to prove that $\displaystyle \exists \frac{d\rho }{dt}\in L^2(0,T;H^{-1}(\mathbb{R}^d))$. Indeed, we have that
$$\frac{\partial \rho }{\partial t}(t,x)=-\nabla \cdot ({\tilde f}(t,x)\rho (t,x))+\frac{\partial ^2(a_{ij}\rho )}{\partial x_i\partial x_j}(t,x)=({\cal A}(t)\rho (t))(x),$$
$\forall t\in (0,T], x\in \mathbb{R}^d$. 
This yields
$$-\nabla \cdot ({\tilde f}\rho)=-\nabla \cdot {\tilde f} \rho -{\tilde f}\cdot \nabla \rho \in L^2(Q_T)\subset L^2(0,T;H^{-1}(\mathbb{R}^d)) \, .$$
On the other hand, we get that there exists a positive constant $\bar{C}$ such that for any $t\in (0,T]$ 
$$\Big| \Big\langle \frac{\partial ^2(a_{ij}(t)\rho (t))}{\partial x_i\partial x_j} ,\varphi \Big\rangle _{H^{-1},H^1}\Big| 
=\Big| \int_{\mathbb{R}^d}\frac{\partial (a_{ij}(t)\rho (t))}{\partial x_j}\frac{\partial \varphi }{\partial x_i}dx\Big| $$
$$\leq \bar{C} \| \rho (t)\| _{H^1}\| \varphi \|_{H^1}, \quad \forall \varphi \in H^1(\mathbb{R}^d) \, .$$
It follows that
$$\frac{\partial ^2(a_{ij}\rho )}{\partial x_i\partial x_j} \in L^2(0,T; H^{-1}(\mathbb{R}^d)) \, .$$
Since ${\cal A}\rho \in L^2(0,T; H^{-1}(\mathbb{R}^d))$ we may infer that for any $\varphi \in H^1(\mathbb{R}^d)$ and $t\in (0,T]$:
$$\langle \rho (t),\varphi \rangle _{L^2}-\langle \rho _0,\varphi \rangle _{L^2}=\langle \rho (t)-\rho _0,\varphi \rangle_{H^{-1},H^1}=\int_0^t\langle {\cal A}(s)\rho (s),\varphi \rangle_{H^{-1}, H^1}ds $$
and so
$$\rho (t)-\rho _0=\int_0^t{\cal A}(s)\rho (s)ds \quad \mbox{\rm in } H^{-1}(\mathbb{R}^d) \, . $$
It follows that $\exists \displaystyle \frac{d\rho }{dt}\in L^2(0,T; H^{-1}(\mathbb{R}^d))$ and that 
$$\frac{d\rho }{dt}(t)={\cal A}(t)\rho(t) \quad \mbox{\rm in } H^{-1}(\mathbb{R}^d), \ \mbox{\rm a.e. } t\in (0,T) \, .$$
It follows that $\rho \in W^{1,2}([0,T]; H^{-1}(\mathbb{R}^d))$ and that $\rho $ is a weak solution to (\ref{eqLFP}). We recall that the weak solution to (\ref{eqLFP}) is unique.
\vspace{3mm}

Let us note that for any $\varepsilon \in (0,T)$, the following function
$$\tilde{\rho }_{\varepsilon }(t,x)
=\int_{\mathbb{R}^d}\Gamma (t,x,\varepsilon ,\xi )\rho (\varepsilon ,\xi )d\xi +\int_{\varepsilon }^t\int_{\mathbb{R}^d}\Gamma (t,x,s,\xi ){\tilde F}(s,\xi )d\xi \ ds $$
$$=\int_{\mathbb{R}^d}\Gamma (t,x,\varepsilon ,\xi )\rho (\varepsilon ,\xi )d\xi +\int_{\varepsilon }^t\int_{\mathbb{R}^d}\nabla _{\xi }\Gamma (t,x,s,\xi )\cdot {\tilde f}(s,\xi )\rho (s,\xi )d\xi \ ds, \quad (t,x)\in (\varepsilon ,T]\times \mathbb{R}^d$$
is a classical solution (see Theorem 4.6, p. 142 in \cite{friedman}) to
\begin{equation}\label{eq1'}
\left\{ \begin{array}{ll}
\displaystyle
\frac{\partial \tilde{\rho }}{\partial t}(t,x)-\frac{\partial ^2(a_{ij}\tilde{\rho })}{\partial x_i\partial x_j}(t,x)=-\nabla \cdot ({\tilde f}(t,x)\rho(t,x)), \quad &(t,x)\in (\varepsilon ,T]\times \mathbb{R}^d , \vspace{2mm} \\
\tilde{\rho }(\varepsilon ,x)=\rho (\varepsilon ,x), & x\in \mathbb{R}^d \, .
\end{array}
\right.
\end{equation}
It follows in the same manner as for $\rho $ that $\tilde{\rho }_{\varepsilon }$ is a weak solution to (\ref{eq1'}).
Since $\rho $ and $\tilde{\rho }_{\varepsilon }$ are both weak solutions to (\ref{eq1'}) and by the uniqueness of this problem we get that $\tilde{\rho }_{\varepsilon }=\rho $ on $[\varepsilon ,T]\times \mathbb{R}^d$
and so
$$\rho (t,x)
=\int_{\mathbb{R}^d}\Gamma (t,x,\varepsilon ,\xi )\rho (\varepsilon ,\xi )d\xi  
+\int_{\varepsilon }^t\int_{\mathbb{R}^d}\nabla _{\xi }\Gamma (t,x,s,\xi )\cdot \tilde{f}(s,\xi )\rho (s,\xi )d\xi \ ds, \quad (t,x)\in (\varepsilon ,T]\times \mathbb{R}^d,$$
for any $\varepsilon \in (0,T)$. Passing to the limit ($\varepsilon \rightarrow 0$) in $L^2(\mathbb{R}^d)$ we get that (\ref{eq9}) holds for any $t \in (0,T]$ and almost any $x\in \mathbb{R}^d$.
\vspace{3mm}

Let us prove now that if the following less restrictive assumptions hold: $\tilde{f}\in L^{\infty }(Q_T)^d$ and $\rho _0\in L^1(\mathbb{R}^d)\cap L^2(\mathbb{R}^d)$, then $\rho $, the unique weak solution to (\ref{eqLFP}), satisfies (\ref{eq9}).
Consider $\{ f_n\} _{n\in \mathbb{N}^*}\subset C_0^{\infty }(\mathbb{R}\times \mathbb{R}^d)^d$ and $\{ \rho _{0n}\} _{n\in \mathbb{N}^*}\subset C_0(\mathbb{R}^d)$
such that $\| f_n\| _{L^{\infty }}\leq \| \tilde{f}\| _{L^{\infty }}$, $f_n\longrightarrow \tilde{f}$ a.e. in $Q_T$ and $\rho _{0n}\longrightarrow \rho _0$ in $L^2(\mathbb{R}^d)$.

Let $\rho _n$ be the classical and also weak solution to (\ref{eqLFP}) corresponding to $\tilde{f}:=f_n$ and $\rho _0:=\rho _{0n}$. It follows in a standard manner that
$$\rho _n\longrightarrow \rho \quad \mbox{\rm in } C([0,T]; L^2(\mathbb{R}^d)) \ \mbox{\rm and in } L^2(0,T; H^1(\mathbb{R}^d)) \, .$$
This yields
\begin{equation}\label{eq9bis}
	\rho _n(t,x)
=\int_{\mathbb{R}^d}\Gamma (t,x,0,\xi )\rho _{0n}(\xi ) \, d\xi  
+\int_0^t\int_{\mathbb{R}^d}\nabla _{\xi }\Gamma (t,x,s,\xi )\cdot f_n(s,\xi )\rho _n(s,\xi ) \, d\xi \ ds \, ,
\end{equation}
for any $(t,x)\in (0,T]\times \mathbb{R}^d$. We will prove that passing to the limit we get (\ref{eq9}).

Let us show firstly that for any $t\in (0,T]$ we get
\begin{equation}\label{conv1}
\int_{\mathbb{R}^d}\Gamma (t,x,0,\xi )\rho _{0n}(\xi ) \, d\xi 
\longrightarrow \int_{\mathbb{R}^d}\Gamma (t,x,0,\xi )\rho _{0}(\xi ) \, d\xi \quad \mbox{\rm in } L^2(\mathbb{R}^d) \, ,
\end{equation}
which is equivalent to
$$\int_{\mathbb{R}^d}\Big| \int_{\mathbb{R}^d}\Gamma (t,x,0,\xi )(\rho _{0n}(\xi )-\rho_0(\xi ))d\xi \Big| ^2dx \longrightarrow 0 \, .$$  
Indeed, we have that
$$\Big\| \int_{\mathbb{R}^d}\Gamma (t, \cdot ,0,\xi )(\rho _{0n}(\xi )-\rho_0(\xi )) \, d\xi \Big\|_{L^2} \leq  
C\int_{\mathbb{R}^d}\frac{e^{-c\frac{|x |^2}{t}}}{t^{\frac{d}{2}}} \, dx \, \| \rho_{0n}-\rho _0\|_{L^2}$$
(by Young's inequality for convolution products)
$$=C_1\| \rho_{0n}-\rho _0\|_{L^2} \longrightarrow 0, $$
where $C_1=C\displaystyle \int_{\mathbb{R}^d}e^{-c|z|^2}dz $, and so we get (\ref{conv1}).

Let us prove that for any $t\in (0,T]$:
\begin{equation}\label{conv2}
\int_0^t\int_{\mathbb{R}^d}\nabla _{\xi }\Gamma (t,x,s,\xi )\cdot f_n(s,\xi )\rho _n(s,\xi ) \, d\xi \, ds \longrightarrow 
\int_0^t\int_{\mathbb{R}^d}\nabla _{\xi }\Gamma (t,x,s,\xi )\cdot \tilde{f}(s,\xi )\rho (s,\xi ) \, d\xi \, ds \ 
\end{equation}
in $L^2(\mathbb{R}^d)$, which is equivalent to
$$\Big\|\int_0^t\int_{\mathbb{R}^d}\nabla _{\xi }\Gamma (t,\cdot ,s,\xi )\cdot (f_n(s,\xi )\rho _n(s,\xi )-\tilde{f}(s,\xi )\rho (s,\xi )) \, d\xi \, ds \Big\|_{L^2} \longrightarrow 0 \, .$$
Indeed, we have that
\begin{equation*}
	\begin{array}{ll}
		\displaystyle
\int_{\mathbb{R}^d}\Big|\int_0^t\int_{\mathbb{R}^d}\nabla _{\xi }\Gamma (t,x,s,\xi )\cdot (f_n(s,\xi )\rho _n(s,\xi )-\tilde{f}(s,\xi )\rho (s,\xi )) \, d\xi \, ds\Big|^2 dx \vspace{2mm} \\
\displaystyle
\leq C^2\int_{\mathbb{R}^d}\int_0^t\int_{\mathbb{R}^d}\frac{e^{-c\frac{|x-\xi |^2}{t-s}}}{(t-s)^{\frac{d+1}{2}}} \, d\xi \, ds 
\int_0^t\int_{\mathbb{R}^d}\frac{e^{-c\frac{|x-\xi |^2}{t-s}}}{(t-s)^{\frac{d+1}{2}}} |f_n(s,\xi )\rho _n(s,\xi )-\tilde{f}(s,\xi )\rho (s,\xi )|^2 \, d\xi \, ds \, dx \vspace{2mm} \\
\displaystyle
\leq C \, C_1\int_{\mathbb{R}^d}\int_0^t\frac{1}{(t-s)^{\frac{1}{2}}}ds 
\int_0^t\int_{\mathbb{R}^d}\frac{e^{-c\frac{|x-\xi |^2}{t-s}}}{(t-s)^{\frac{d+1}{2}}} |f_n(s,\xi )\rho _n(s,\xi )-\tilde{f}(s,\xi )\rho (s,\xi )|^2 \, d\xi \, ds \, dx \vspace{2mm} \\
\displaystyle
=2C \, C_1t^{\frac{1}{2}}\int_{\mathbb{R}^d}\int_0^t
|f_n(s,\xi )\rho _n(s,\xi )-\tilde{f}(s,\xi )\rho (s,\xi )|^2
\int_{\mathbb{R}^d}
\frac{e^{-c\frac{|x-\xi |^2}{t-s}}}{(t-s)^{\frac{d+1}{2}}} \, dx \, ds \, d\xi \vspace{2mm} \\
\displaystyle
=2C_1^2t^{\frac{1}{2}}\int_0^t\frac{1}{(t-s)^{\frac{1}{2}}}
\|f_n(s)\rho _n(s)-\tilde{f}(s)\rho (s)\|^2_{L^2} \, ds \, .
\end{array}
\end{equation*}
On the other hand, for almost any $s\in (0,t)$: 
$$f_n(s)\rho _n(s)-\tilde{f}(s) \rho (s)=f_n(s)(\rho_n(s)-\rho (s))+(f_n(s)-\tilde{f}(s))\rho (s)\longrightarrow 0 \ \mbox{\rm in } L^2(\mathbb{R}^d)^d$$
(because $\{ f_n(s)\}_{n\in \mathbb{N}^*}$ is bounded in $L^{\infty }(\mathbb{R}^d)^d$ and $\{ \rho _n(s)- \rho (s)\}_{n\in \mathbb{N}^*}$ converges to $0$ in $L^2(\mathbb{R}^d)$ and on the other hand
$\{ (f_n(s)-\tilde{f}(s))\rho (s)\}_{n\in \mathbb{N}^*} $ converges to $0$ in $L^2(\mathbb{R}^d)^d$ via Lebesgue's dominated convergence theorem). So,
$$\frac{1}{(t-s)^{\frac{1}{2}}}
\|f_n(s)\rho _n(s)-\tilde{f}(s)\rho (s)\|^2_{L^2}\longrightarrow 0$$
a.e. $s\in (0,t)$.
This yields
$$\| f_n(s)\rho_n(s)-\tilde{f}(s)\rho (s)\| _{L^2}^2\leq 2\| f_n(s)(\rho _n(s)-\rho (s))\|^2_{L^2}
+2\| (f_n(s)-\tilde{f}(s))\rho (s)\| ^2_{L^2}$$
$$\leq 2\| \tilde{f}\|^2_{L^{\infty }}\| \rho _n-\rho \|^2_{C([0,T];L^2(\mathbb{R}^d))}
+8\| \tilde{f}\|_{L^{\infty }}^2\| \rho \| ^2_{C([0,T];L^2(\mathbb{R}^d))}\leq M_1,$$
where $M_1$ is a positive constant (because $\{ \rho _n-\rho \}_{n\in \mathbb{N}^*}$ is bounded in $C([0,T];L^2(\mathbb{R}^d))$).
We infer that
$$\frac{1}{(t-s)^{\frac{1}{2}}}\|f_n(s)\rho _n(s)-\tilde{f}(s)\rho (s)\|^2_{L^2}
\leq \frac{M_1}{(t-s)^{\frac{1}{2}}}$$
a.e. (and the function in the right-hand side belongs to $L^1(0,t)$) and so we may conclude via Lebesgue's dominated convergence theorem that
$$\int_0^t\frac{1}{(t-s)^{\frac{1}{2}}}
\|f_n(s)\rho _n(s)-\tilde{f}(s)\rho (s)\|^2_{L^2}ds \longrightarrow 0$$
and consequently (\ref{conv2}) holds. 

If we pass to the limit in (\ref{eq9bis}) (in $L^2(\mathbb{R}^d)$) we conclude that (\ref{eq9}) holds for and $t\in (0,T]$ in $L^2(\mathbb{R}^d)$ and also a.e. 
$(t,x)\in Q_T$.

\vspace{5mm}

Assume now that additionally $\rho _0\in L^m(\mathbb{R}^d)$, for $m\in [1, +\infty )$ arbitrary but fixed. Consider the operator
$K:L^m(Q_T)\longrightarrow L^m(Q_T)$, given by
\begin{equation*}
K(y)(t,x)=\int_{\mathbb{R}^d}\Gamma (t,x,0,\xi )\rho _0(\xi ) \, d\xi +\int_0^t\int_{\mathbb{R}^d}\nabla _{\xi }\Gamma (t,x,s,\xi )\cdot \tilde{f}(s,\xi )y(s,\xi ) \, d\xi \ ds ,
\end{equation*} 
for $y\in L^m(Q_T)$. The properties of $\Gamma $ allows us to conclude that $K$ is correctly defined.

We will show that $K$ is a contraction if we consider on $L^m(Q_T)$ the following norm (equivalent to the usual one)
$$\| | y| \|=\Big( \int_{Q_T}e^{-\lambda t}|y(t,x)|^mdx \, dt\Big)^{\frac{1}{m}} ,$$
where $\lambda $ is a positive constant to be precised later.
\vspace{2mm}

Assume for the beginning $m\in (1,+\infty )$. 
Let $m^*$ such that $\frac{1}{m}+\frac{1}{m^*}=1$.
For any $y\in L^m(Q_T)$ we have that
\begin{equation*}
	\begin{array}{ll}
\displaystyle 
\int_0^T\int_{\mathbb{R}^d}e^{-\lambda t}\Big| \int_0^t\int_{\mathbb{R}^d}\nabla _{\xi }\Gamma (t,x,s,\xi )\cdot \tilde{f}(s,\xi )y(s,\xi )d\xi \, ds \Big| ^m \, dx \, dt 
\vspace{2mm} \\
\displaystyle  
\leq C^m\| \tilde{f}\|_{L^{\infty }}^m\int_0^T\int_{\mathbb{R}^d}e^{-\lambda t}\Big| \int_0^t\int_{\mathbb{R}^d}\frac{e^{-c \frac{|x-\xi|^2}{t-s}}}{(t-s)^{\frac{d+1}{2}}}|y(s,\xi )|d\xi \, ds \Big| ^m \, dx \, dt \vspace{2mm} \\
\displaystyle \leq \| \tilde{f}\|_{L^{\infty }}^m C^m\int_0^T\int_{\mathbb{R}^d}e^{-\lambda t}\int_0^t\int_{\mathbb{R}^d}\frac{e^{-c \frac{|x-\xi|^2}{t-s}}}{(t-s)^{\frac{d+1}{2}}}|y(s,\xi )|^md\xi \, ds \Big( \int_0^t\int_{\mathbb{R}^d}\frac{e^{-c \frac{|x-\xi|^2}{t-s}}}{(t-s)^{\frac{d+1}{2}}}d\xi  \, ds \Big) ^{\frac{m}{m^*}}dx \, dt 
\vspace{2mm} \\
\displaystyle \leq \| \tilde{f}\|_{L^{\infty }}^mC \, C_12^{\frac{m}{m^*}}t^{\frac{m}{2m^*}}\int_0^T\int_{\mathbb{R}^d}e^{-\lambda t}\int_0^t\int_{\mathbb{R}^d}\frac{1}{(t-s)^{\frac{d+1}{2}}}e^{-c \frac{|x-\xi|^2}{t-s}}|y(s,\xi )|^md\xi \, ds \, dx \, dt \vspace{2mm} \\
\displaystyle 
\leq \| \tilde{f}\|_{L^{\infty }}^mC \, C_12^{\frac{m}{m^*}}t^{\frac{m}{2m^*}}\int_0^T\int_{\mathbb{R}^d}e^{-\lambda s}|y(s,\xi )|^m\int_s^T\frac{1}{(t-s)^{\frac{d+1}{2}}}e^{-\lambda (t-s)}\int_{\mathbb{R}^d}e^{-c \frac{|x-\xi|^2}{t-s}}dx \, dt \, d\xi \, ds \vspace{2mm} \\
\displaystyle
\leq \| \tilde{f}\|_{L^{\infty }}^mC_1^22^{\frac{m}{m^*}}t^{\frac{m}{2m^*}}\int_0^T\int_{\mathbb{R}^d}e^{-\lambda s}|y(s,\xi )|^m\int_s^T\frac{e^{-\lambda (t-s)}}{(t-s)^{\frac{1}{2}}} dt \, d\xi \, ds \, .
\end{array}
\end{equation*}
For any $t\in (s,T]$ we get that
$$e^{\lambda (t-s)}\geq (4\lambda )^{\frac{1}{4}}(t-s)^{\frac{1}{4}}\Longrightarrow \frac{e^{-\lambda (t-s)}}{(t-s)^{\frac{1}{2}}}\leq \frac{1}{(4\lambda )^{\frac{1}{4}}(t-s)^{\frac{3}{4}}}$$
and so there exists a positive constant $C_2$ such that
$$\int_0^T\int_{\mathbb{R}^d}e^{-\lambda t}\Big| \int_0^t\int_{\mathbb{R}^d}\nabla _{\xi }\Gamma (t,x,s,\xi )\cdot \tilde{f}(s,\xi )y(s,\xi )d\xi \, ds \Big| ^mdx \, dt \leq \frac{C_2}{\lambda ^{\frac{1}{4}}}\int_0^T\int_{\mathbb{R}^d}e^{-\lambda s}|y(s,\xi )|^md\xi \, ds.$$
We may infer that for any $y_1, y_2\in L^m(Q_T)$ we get that
$$\| | K(y_1-y_2)| \| \leq \Big( \frac{C_2}{\lambda ^{\frac{1}{4}}}\Big)^{\frac{1}{m}}\| | y_1-y_2| \| ,$$
where $C_2$ is a positive constant independent of $y_1, y_2$ and $\lambda $. 

It follows that if we consider $\lambda $ a constant greater that $C_2^4$, then $K$ is a contraction. We get via Banach's fixed point theorem that $K$ admits a unique fixed point $\tilde{\rho }\in L^m(Q_T)$. 
Since the fixed point of $K$ and the weak solution to (\ref{eqLFP}) are limits (in different functional spaces) of the same sequence of iterations we conclude that $\tilde{\rho }(t,x)=\rho (t,x)$ a.e. $(t,x)\in Q_T$.

We get that there exists positive constants $c_9, c_{10}, c_{11}, c_{12}, c_{13}, c_{14}$ such that for any $t\in [0,T]$:
\begin{equation*}
	\begin{array}{ll}
\| \rho (t)\|_{L^m}^m &\displaystyle \leq c_9\Big[\int_{\mathbb{R}^d}\Big| \int_{\mathbb{R}^d}\Gamma (t,x,0,\xi )\rho _0(\xi )d\xi \Big| ^mdx \vspace{2mm} \\
 ~ & ~~~\displaystyle +\int_{\mathbb{R}^d}\Big| \int_0^t\int_{\mathbb{R}^d}\nabla _{\xi }\Gamma (t,x,s,\xi )\cdot \tilde{f}(s,\xi )\rho (s,\xi ) \, d\xi \ ds \Big|^m dx \Big] \vspace{2mm} \\ 
~ &\leq \displaystyle c_{10}\Big[\int_{\mathbb{R}^d}\Big| \int_{\mathbb{R}^d}\frac{e^{-c\frac{|x-\xi |^2}{t}}}{t^{\frac{d}{2}}}\rho _0(\xi ) \, d\xi \Big| ^mdx
+\int_{\mathbb{R}^d}\Big| \int_0^t\int_{\mathbb{R}^d}\frac{e^{-c\frac{|x-\xi |^2}{t-s}}}{(t-s)^{\frac{d+1}{2}}}\rho (s,\xi ) \, d\xi \ ds \Big|^m dx \Big] \vspace{2mm} \\
~ & \displaystyle \leq c_{11}\Big[\int_{\mathbb{R}^d}\int_{\mathbb{R}^d}\frac{e^{-c\frac{|x-\xi |^2}{t}}}{t^{\frac{d}{2}}}|\rho _0(\xi )|^md\xi \, \Big(\int_{\mathbb{R}^d}\frac{e^{-c\frac{|x-\xi |^2}{t}}}{t^{\frac{d}{2}}} d\xi \Big)^{\frac{m}{m^*}} \, dx \vspace{2mm} \\
~ & ~~~\displaystyle +\int_{\mathbb{R}^d} \int_0^t\int_{\mathbb{R}^d}\frac{e^{-c\frac{|x-\xi |^2}{t-s}}}{(t-s)^{\frac{d+1}{2}}}|\rho (s,\xi )|^md\xi \ ds  \, \Big(\int_0^t\int_{\mathbb{R}^d}\frac{e^{-c\frac{|x-\xi |^2}{t-s}}}{(t-s)^{\frac{d+1}{2}}}d\xi \ ds \Big)^{\frac{m}{m^*}} \, dx \Big] \vspace{2mm} \\
~ & \displaystyle \leq c_{12}\Big[\int_{\mathbb{R}^d}|\rho _0(\xi )|^m \, \int_{\mathbb{R}^d}\frac{e^{-c\frac{|x-\xi |^2}{t}}}{t^{\frac{d}{2}}}dx  \, d\xi +\int_{\mathbb{R}^d} \int_0^t\int_{\mathbb{R}^d}|\rho (s,\xi )|^m \,  \frac{e^{-c\frac{|x-\xi |^2}{t-s}}}{(t-s)^{\frac{d+1}{2}}}dx \ ds \ d\xi \Big] \vspace{2mm} \\
~ & \displaystyle \leq c_{13}\Big[\|\rho _0\| _{L^m}^m +\int_0^t\frac{\| \rho (s)\| _{L^m}^m}{(t-s)^{\frac{1}{2}}} ds  \Big] \, \vspace{2mm} \\
~ & \displaystyle \leq c_{13}\Big[\|\rho _0\| _{L^m}^m +\Big( \int_0^t\| \rho (s)\| _{L^m}^{3m} ds\Big) ^{\frac{1}{3}}\Big(\int_0^t\frac{1}{(t-s)^{\frac{1}{2}\frac{3}{2}}} ds\Big)^{\frac{2}{3}}  \Big] \vspace{2mm} \\
~ & \displaystyle \leq c_{14}\Big[\|\rho _0\| _{L^m}^m +\Big(\int_0^t\| \rho (s)\| _{L^m}^{3m} ds \Big)^{\frac{1}{3}} \Big] \, ,
\end{array}
\end{equation*}
where $\rho $ is the weak solution to (\ref{eqLFP}).
It follows that there exists a positive constant $c_{15}$ such that for any $t\in [0,T]$:
$$\| \rho (t)\|_{L^m}^{3m}\leq c_{15}\Big[\|\rho _0\| _{L^m}^{3m} +\int_0^t\| \rho (s)\| _{L^m}^{3m} ds  \Big] \, $$
and by Gronwall's inequality that
$$\| \rho (t)\|_{L^m}^{3m}\leq c_{16}\| \rho _0\|_{L^m}^{3m} \, ,$$
where $c_{16}$ is a positive constant and consequently there exists a positive constant $c_T$ such that for any $t\in [0,T]$:
$$\| \rho (t)\|_{L^m}\leq c_T\| \rho _0\|_{L^m}\, .$$
\vspace{3mm}

Assume now that $m=1$. We adapt the past argument and get that for any $y\in L^1(Q_T)$ we have that
\begin{equation*}
	\begin{array}{ll}
\displaystyle 
\int_0^T\int_{\mathbb{R}^d}e^{-\lambda t}\Big| \int_0^t\int_{\mathbb{R}^d}\nabla _{\xi }\Gamma (t,x,s,\xi )\cdot \tilde{f}(s,\xi )y(s,\xi ) \, d\xi \, ds \Big| dx \, dt 
\vspace{2mm} \\
\displaystyle  
\leq C\| \tilde{f}\|_{L^{\infty }}\int_0^T\int_{\mathbb{R}^d}e^{-\lambda t}\Big| \int_0^t\int_{\mathbb{R}^d}\frac{1}{(t-s)^{\frac{d+1}{2}}}e^{-c \frac{|x-\xi|^2}{t-s}}|y(s,\xi )| \, d\xi \, ds \Big| dx \, dt \vspace{2mm} \\
\displaystyle 
\leq C\| \tilde{f}\|_{L^{\infty }}\int_0^T\int_{\mathbb{R}^d}e^{-\lambda s}|y(s,\xi )|\int_s^T\frac{1}{(t-s)^{\frac{d+1}{2}}}e^{-\lambda (t-s)}\int_{\mathbb{R}^d}e^{-c \frac{|x-\xi|^2}{t-s}}dx \, dt \, d\xi \, ds \vspace{2mm} \\
\displaystyle
\leq C_1\| \tilde{f}\|_{L^{\infty }}\int_0^T\int_{\mathbb{R}^d}e^{-\lambda s}|y(s,\xi )|\int_s^T\frac{e^{-\lambda (t-s)}}{(t-s)^{\frac{1}{2}}} dt \, d\xi \, ds \, .
\end{array}
\end{equation*}
For any $t\in (s,T]$ we get that $\frac{e^{-\lambda (t-s)}}{(t-s)^{\frac{1}{2}}}\leq \frac{1}{(4\lambda )^{\frac{1}{4}}(t-s)^{\frac{3}{4}}}$
and so there exists a positive constant $C_3$ such that
$$\int_0^T\int_{\mathbb{R}^d}e^{-\lambda t}\Big| \int_0^t\int_{\mathbb{R}^d}\nabla _{\xi }\Gamma (t,x,s,\xi )\cdot \tilde{f}(s,\xi )y(s,\xi ) \, d\xi \, ds \Big| dx \, dt \leq \frac{C_3}{\lambda ^{\frac{1}{4}}}\int_0^T\int_{\mathbb{R}^d}e^{-\lambda s}|y(s,\xi )|d\xi \, ds.$$
We may infer that for any $y_1, y_2\in L^1(Q_T)$ we get that
$$\| | K(y_1-y_2)| \| \leq  \frac{C_3}{\lambda ^{\frac{1}{4}}}\| | y_1-y_2| \| ,$$
where $C_3$ is a positive constant independent of $y_1, y_2$ and $\lambda $. 

It follows that if we consider $\lambda $ a constant greater that $C_3^4$, then $K$ is a contraction. We get via Banach's fixed point theorem that $K$ admits a unique fixed point $\tilde{\rho }\in L^1(Q_T)$. 
Since the fixed point of $K$ and the weak solution to (\ref{eqLFP}) are limits (in different functional spaces) of the same sequence of iterations we conclude that $\tilde{\rho }(t,x)=\rho (t,x)$ a.e. $(t,x)\in Q_T$.

We get that there exist the positive constants $c_{17}, c_{18}, c_{19}$ such that for any $t\in [0,T]$:
\begin{equation*}
	\begin{array}{ll}
\| \rho (t)\|_{L^1} &\displaystyle \leq \int_{\mathbb{R}^d}\Big| \int_{\mathbb{R}^d}\Gamma (t,x,0,\xi )\rho _0(\xi ) \, d\xi \Big| dx \vspace{2mm} \\
 ~ & ~~~\displaystyle +\int_{\mathbb{R}^d}\Big| \int_0^t\int_{\mathbb{R}^d}\nabla _{\xi }\Gamma (t,x,s,\xi )\cdot \tilde{f}(s,\xi )\rho (s,\xi ) \, d\xi \ ds \Big| dx  \vspace{2mm} \\ 
~ &\leq \displaystyle c_{17}\Big[\int_{\mathbb{R}^d}\Big| \int_{\mathbb{R}^d}\frac{e^{-c\frac{|x-\xi |^2}{t}}}{t^{\frac{d}{2}}}\rho _0(\xi ) \, d\xi \Big| dx
+\int_{\mathbb{R}^d}\Big| \int_0^t\int_{\mathbb{R}^d}\frac{e^{-c\frac{|x-\xi |^2}{t-s}}}{(t-s)^{\frac{d+1}{2}}}\rho (s,\xi ) \, d\xi \ ds \Big| dx \Big] \vspace{2mm} \\
~ & \displaystyle \leq c_{17}\Big[\int_{\mathbb{R}^d}|\rho _0(\xi )| \, \int_{\mathbb{R}^d}\frac{e^{-c\frac{|x-\xi |^2}{t}}}{t^{\frac{d}{2}}}dx  \, d\xi +\int_{\mathbb{R}^d} \int_0^t\int_{\mathbb{R}^d}|\rho (s,\xi )| \,  \frac{e^{-c\frac{|x-\xi |^2}{t-s}}}{(t-s)^{\frac{d+1}{2}}}dx \ ds \ d\xi \Big] \vspace{2mm} \\
~ & \displaystyle \leq c_{18}\Big[\|\rho _0\| _{L^1} +\int_0^t\frac{\| \rho (s)\| _{L^1}}{(t-s)^{\frac{1}{2}}} ds  \Big] \, \vspace{2mm} \\
~ & \displaystyle \leq c_{18}\Big[\|\rho _0\| _{L^1} +\Big( \int_0^t\| \rho (s)\| _{L^1}^{3} ds\Big) ^{\frac{1}{3}}\Big(\int_0^t\frac{1}{(t-s)^{\frac{1}{2}\frac{3}{2}}} ds\Big)^{\frac{2}{3}}  \Big] \vspace{2mm} \\
~ & \displaystyle \leq c_{19}\Big[\|\rho _0\| _{L^1}+\Big(\int_0^t\| \rho (s)\| _{L^1}^{3} ds \Big)^{\frac{1}{3}} \Big] \, ,
\end{array}
\end{equation*}
where $\rho $ is the weak solution to (\ref{eqLFP}).
It follows that there exists a positive constant $c_{20}$ such that for any $t\in [0,T]$:
$$\| \rho (t)\|_{L^1}^{3}\leq c_{20}\Big[\|\rho _0\| _{L^1}^{3} +\int_0^t\| \rho (s)\| _{L^1}^{3} ds  \Big] \, $$
and by Gronwall's inequality that
$$\| \rho (t)\|_{L^1}^{3}\leq c_{21}\| \rho _0\|_{L^1}^{3} \, ,$$
where $c_{21}$ is a positive constant and consequently there exists a positive constant $c_T$ such that for any $t\in [0,T]$:
$$\| \rho (t)\|_{L^1}\leq c_T\| \rho _0\|_{L^1}\, .$$
\vspace{3mm}

Assume now that additionally $\rho _0\in L^{\infty }(\mathbb{R}^d)$. Consider the operator
$K:L^{\infty }(Q_T)\longrightarrow L^{\infty }(Q_T)$, given by
\begin{equation*}
K(y)(t,x)=\int_{\mathbb{R}^d}\Gamma (t,x,0,\xi )\rho _0(\xi ) \, d\xi +\int_0^t\int_{\mathbb{R}^d}\nabla _{\xi }\Gamma (t,x,s,\xi )\cdot \tilde{f}(s,\xi )y(s,\xi ) \, d\xi \ ds ,
\end{equation*} 
for $y\in L^{\infty }(Q_T)$. The properties of $\Gamma $ allows us to conclude that $K$ is correctly defined.

It is possible to show that $K$ is a contraction if we consider on $L^{\infty }(Q_T)$ the following norm (equivalent to the usual one)
$$\| | y| \|=Ess ~sup_{(t,x)\in Q_T}e^{-\lambda t}|y(t,x)|,$$
where $\lambda $ is a sufficiently large positive constant.
In a similar manner as for $m\in [1,+\infty )$ it follows that there exists a positive constant $c_T$ such that for any $t\in [0,T]$:
$$\| \rho (t)\|_{L^{\infty }}\leq c_T\| \rho _0\|_{L^{\infty }}\, .$$
\vspace{3mm}

\section{Existence of a solution to the MFG system}

Assume that (H1), (H2) hold.

Let us denote 
$$J(u,\rho )=\int_{Q_T}[L(t,x,u(t,x))\rho (t,x)+G(t,x,\rho (t,x)) ]dx \, dt +\int_{\mathbb{R}^d}G^0(x,\rho (T,x))dx ,$$
$u\in L^{\infty }(Q_T)^{\ell }, \rho \in L^2(Q_T) .$
We may rewrite the optimal control problem as
$$Minimize \, \{ J(u,\rho ^u); \ u \in L^{\infty }(Q_T)^{\ell }\} \, ,$$
\vspace{3mm}

where we assume that $L, g, g_0$ are measurable functions such that

\begin{itemize}
\item[\bf (H3)] $L:Q_T\times \mathbb{R}^{\ell }\longrightarrow \mathbb{R}$ is convex and lower semicontinuous in $u$ for almost any $(t,x)\in Q_T$, measurable in $(t,x)\in Q_T$ for any $u\in \mathbb{R}^{\ell}, |u|\leq a$, and
\begin{equation*}
L(t,x,u)=+\infty \quad \mbox{\rm a.e. } (t,x)\in Q_T, \forall u\in \mathbb{R}^{\ell }, \ |u|>a,
\end{equation*}
where $a$ is a positive constant,
\begin{equation*}
\sup \ \{ |L(t,x,u)|; \ |u|\leq a\} \in L^2(Q_T), \ L(\cdot , \cdot , 0)\in L^{\infty }(Q_T) \, .
\end{equation*}
The assumptions on $L$ imply that $U_0=\{ v\in \mathbb{R}^{\ell }; \ |v|\leq a\}$.

\item[\bf (H4)] $G:Q_T\times \mathbb{R}\longrightarrow \mathbb{R}$ is convex and continuous in $r\in \mathbb{R}$ for almost any $(t,x)\in Q_T$, 
measurable in $(t,x)\in Q_T$ for any $r\in \mathbb{R}$, and 
$$ \alpha (t,x)r\leq G(t,x,r)\leq \bar{C}_1r^2+\bar{\alpha }(t,x), \quad \mbox{\rm a.e. } (t,x)\in Q_T, \ \forall r\in \mathbb{R} ,$$
where $\alpha , \bar{\alpha }\in L^1(Q_T)\cap L^{\infty }(Q_T)$ and $\bar{C}_1$ is a positive constant. 
\item[\bf(H5)] $G^0:\mathbb{R}^d\times \mathbb{R}\longrightarrow \mathbb{R}$ is convex and continuous in $r\in \mathbb{R}$ for almost any $x\in \mathbb{R}^d$, 
measurable in $x\in \mathbb{R}^d$ for any $r\in \mathbb{R}$, and 
$$\alpha _0(x)r\leq G^0(x,r)\leq \bar{C}_{10}r^2+\bar{\alpha }_0(x) \quad \mbox{\rm a.e. } x\in \mathbb{R}^d, \ \forall r\in \mathbb{R} ,$$
where $\alpha _0\in L^1(\mathbb{R}^d)\cap L^{\infty }(Q_T), \ \bar{\alpha }_0\in L^1(Q_T)$ and $\bar{C}_{10}$ is a positive constant.
\vspace{2mm}

$f: [0,T]\times \mathbb{R}^d\times \mathbb{R}^{\ell }$ is a Borel function and satisfies

\item[\bf (H6)] $f(t,x,u)=f_0(t,x)u$, a.e. $(t,x,u)\in {\cal K}=[0,T]\times \mathbb{R}^d\times \left\{ v\in \mathbb{R}^{\ell }; \, |v|\leq a\right\} ,$ and 
$f_0\in L^{\infty }(Q_T)^{d\times {\ell}}$

or 

\item[\bf (H6)'] For almost any $(t,x)\in Q_T$, $f(t,x, \cdot ) \in C^1(\{v\in \mathbb{R}^{\ell}; \ |v|\leq a\}; \mathbb{R}^d)$ and $f$ and $f_u$ are bounded on ${\cal K}
=[0,T]\times \mathbb{R}^d\times \left\{ v\in \mathbb{R}^{\ell }; \, |v|\leq a\right\}$.
The set
$${\cal V}=\Big\{v; \  \exists u\in {\cal U}: \ v(t,x)=f(t,x,u(t,x)), \ \mbox{\rm a.e. } (t,x)\in Q_T\Big\} \, $$
is weak * closed in $L^{\infty }(Q_T)^d$

\noindent
(note that $f$ satisfying (H6) also verifies these assumtions in (H6)').

One moreover assumes that for almost any $(t,x)\in Q_T$, the set $V_{t,x}=\{ f(t,x,u); \ |u|\leq a\}$ is convex and the mapping $u\mapsto f(t,x,u)$ is a homeomorphism between $\{ u\in \mathbb{R}^{\ell}; \ |u|\leq a\}$ and 
$V_{t,x}$ (by abuse of notation we consider that $f^{-1}(t,x,v)$ is the inverse of $u\mapsto f(t,x,u)$) and for almost any $(t,x)\in Q_T$, $L(t,x,f^{-1}(t,x,v))$ is convex with respect to $v\in V_{t,x}$. 
\end{itemize}
\vspace{2mm}

\noindent
\bf Remark. \rm Function $f$ defined by $f(t,x,u)=|u|u $, for $u\in \mathbb{R}^d$ (here $d={\ell }$) satisfies (H6)' if one assumes the convexity of $L(t,x,|v|^{-\frac{1}{2}}v)$ with respect to $v$. 
\vspace{3mm}

Let us recall first that the convex conjugate of $v\mapsto L(t,x,v)$, 
$H:Q_T\times \mathbb{R}^{\ell }\longrightarrow \mathbb{R}$ is convex in $q\in \mathbb{R}^{\ell }$ for almost any $(t,x)\in Q_T$,
measurable in $(t,x)\in Q_T$ for any $q\in \mathbb{R}^{\ell}$ and let us prove that $H$ is Lipschitz continuous with respect to $q$ (for almost any $(t,x)\in Q_T$).
Indeed, for any $q, \bar{q}\in \mathbb{R}^{\ell }$:
$$H(t,x,\bar{q})=\sup \ \left\{ q\cdot u-L(t,x,u)+(\bar{q}-q)\cdot u; \ u\in \mathbb{R}^{\ell }, |u|\leq a\right \} $$
$$\leq H(t,x,q)+\sup \ \left\{ (\bar{q}-q)\cdot u; \ |u|\leq a\right\} \leq H(t,x,q)+a|\bar q-q| $$
and so
$$| H(t,x,\bar{q})-H(t,x,q)|\leq a|\bar{q}-q| \, , $$
which means that $q\mapsto H(t,x,q)$ is Lipschitz continuous.
\vspace{5mm}

Denote by
\begin{equation*}
H_q(t,x,q)=\partial _qH(t,x,q)=\Big\{ \zeta \in \mathbb{R}^{\ell }; \zeta \cdot (\bar{q}-q)\leq H(t,x,\bar{q})-H(t,x,q), \ \forall \bar{q}\in \mathbb{R}^{\ell }\Big\} ,
\end{equation*}
for any $q\in \mathbb{R}^{\ell }$, a.e. $(t,x)\in Q_T$, 
the subdifferential with respect to $q$. We have that for any $\zeta \in H_q(t,x,q)$:
$$\zeta \cdot (\bar{q}-q)\leq H(t,x,\bar{q})-H(t,x,q)\leq a|\bar{q}-q|, \quad \forall q, \bar{q}\in \mathbb{R}^{\ell }$$
and so $|\zeta |\leq a$ and
\begin{equation*}
\sup \ \Big\{ |\zeta |; \ \zeta \in H_q(t,x,q)\Big\} \leq a, \quad \forall q\in \mathbb{R}^{\ell }, \ \mbox{\rm a.e. } (t,x)\in Q_T \, .
\end{equation*}
\vspace{3mm}

\noindent
\bf Theorem 3.1. \it Assume that $\rho _0(x)>0$ a.e. $x\in \mathbb{R}^d$,
$\log \rho _0\in L^1_{loc}(\mathbb{R}^d)$.
Then there exists an optimal pair $(u^*, \rho ^*)$ for problem (P), which satisfies
\begin{equation*}
u^*\in L^{\infty }(Q_T)^{\ell }, \quad \rho ^*(t,x)>0, \ \mbox{\rm a.e. } (t,x)\in Q_T,
\end{equation*}
\begin{equation*}
\rho ^*\in L^2(0,T;H^1(\mathbb{R}^d))\cap C([0,T];L^2(\mathbb{R}^d)), \ \frac{d\rho ^*}{dt}\in L^2(0,T; H^{-1}(\mathbb{R}^d)),
\end{equation*}
\begin{equation*}
u^*(t,x)\in H_q(t,x,f_u^{u^*}(t,x)^T\nabla p(t,x)), \quad \mbox{\rm a.e. } (t,x)\in Q_T,
\end{equation*}
where 
\begin{equation*}
p\in L^2(0,T;H^1(\mathbb{R}^d))\cap W^{1,2}([0,T];H^{-1}(\mathbb{R}^d))  
\end{equation*}
is a weak solution to
\begin{equation}\label{eqHJp}
\left\{ \begin{array}{ll}
\displaystyle
\frac{\partial p}{\partial t}(t,x)+a_{ij}(t,x)\frac{\partial ^2p}{\partial x_i\partial x_j}(t,x)+H(t,x,f_u^{u^*}(t,x)^T \nabla p(t,x)) & \vspace{2mm} \\
\hspace{3cm} = (f_u^{u^*}(t,x)u^*(t,x)-f^{u^*}(t,x))\cdot \nabla p(t,x)+\eta (t,x) , \ & (t,x)\in Q_T, \vspace{2mm} \\
p(T,x)=-\eta _0(x), & x\in \mathbb{R}^d \, . 
\end{array}
\right.
\end{equation}
Here 
\begin{equation*}
\eta \in L^2(Q_T), \ \eta (t,x)\in G_r(t,x,\rho ^*(t,x)), \quad \mbox{\rm a.e. } (t,x)\in Q_T , 
\end{equation*}
\begin{equation*}
\eta _0\in L^2(\mathbb{R}^d), \ \eta _0(x)\in G^0_r(x, \rho ^*(T,x)), \quad \mbox{\rm a.e. } x\in \mathbb{R}^d \, .
\end{equation*}
\vspace{2mm}

\it Proof. \rm Let us prove first the existence of an optimal control. We prove first that
$J(u,\rho ^u)$ is not identically $+\infty $. Indeed, for any $u\in {\cal U}=\{ v\in L^{\infty }(Q_T)^{\ell }; \ |v(t,x)| \leq a, \ \mbox{\rm a.e. in } Q_T\}$:
Using (H4) we get that
$$\alpha (t,x)\rho ^u(t,x)\leq G(t,x,\rho ^u(t,x))\leq \bar{C}_1\rho ^u(t,x)^2+\bar{\alpha }(t,x)$$
for $(t,x)\in Q_T$. It follows that
$$G(t,x,\rho ^u(t,x))\in L^1(Q_T) \, .$$

We get that for any $u\in U$: $L(t,x,u(t,x))\in L^2(Q_T)$ and so $L(t,x,u(t,x))\rho ^u(t,x)\in L^1(Q_T)$.
We have that for any $u\in U$:
$$\alpha _0(x)\rho ^u(T,x)\leq G^0(x,\rho ^u(T,x))\leq \bar{C}_{10}\rho ^u(T,x)^2+\bar{\alpha }_0(x)$$
and this implies that $G^0(x,\rho ^u(T,x))\in L^1(\mathbb{R}^d)$.

We conclude that indeed for any $u\in {\cal U}$: $J(u,\rho ^u)\in \mathbb{R}$.
\vspace{2mm}

On the other hand for any $u\in {\cal U}$: 
$$J(u,\rho ^u)\geq -\int_{Q_T}\sup _{|v|\leq a}|L(t,x,v)| \rho ^u(t,x)dx \ dt +\int_{Q_T}\alpha (t,x)\rho ^u(t,x)dx \ dt$$
$$+\int_{\mathbb{R}^d}\alpha _0(x)\rho ^u(T,x)dx\geq -C^0,$$
where $C^0$ is a real constant.
It follows that 
$$\inf _{u\in U}J(u,\rho ^u)\in \mathbb{R} \, .$$
\vspace{5mm}

We shall use the notation $\rho^{0,v}$ for the weak solution to 
\begin{equation*}
\left\{ \begin{array}{ll}
\displaystyle
\frac{\partial \rho}{\partial t}(t,x)=-\nabla \cdot (v(t,x) \, \rho(t,x))+\frac{\partial ^2(a_{ij}\rho )}{\partial x_i\partial x_j}(t,x), \quad &(t,x)\in Q_T,\vspace{2mm} \\
\rho(0,x)=\rho_0(x), & x\in \mathbb{R}^d \, .
\end{array}
\right.
\end{equation*}

\noindent
Let $\{ (u_n,\rho ^{u_n})\}_{n\in \mathbb{N}^*}$ such that
$$\inf _{u\in U}J(u,\rho ^u)\leq J(u_n,\rho ^{u_n})<\inf _{u\in U}J(u,\rho ^u)+\frac{1}{n}, \quad \forall n\in \mathbb{N}^* \, .$$
Denote by $v_n(t,x)=f(t,x,u_n(t,x))$ for $(t,x)\in Q_T$. So, we get
$$\rho ^{0,v_n}=\rho ^{u_n} \, .$$

Using Theorem 2.1 we get that on a subsequence (also indexed by $n$):
$$v_n\longrightarrow v^* \quad \mbox{\rm weak * in } L^{\infty }(Q_T)^d , $$
$$\rho ^{0,v_n}\longrightarrow \rho ^{0,v^*} \quad \mbox{\rm weakly in } L^2(0,T; H^1(\mathbb{R}^d)), $$
and via Aubin's compactness theorem, also strongly in $L^2(0,T; L^2_{loc}(\mathbb{R}^d))$, on a subsequence (also indexed by $n$)
$$\rho ^{0,v_n}(T)\longrightarrow \rho ^{0,v^*}(T) \quad \mbox{\rm weakly in } L^2(\mathbb{R}^d) ,$$
$$\frac{d\rho ^{0,v_n}}{dt}\longrightarrow \frac{d\rho ^{0,v^*}}{dt} \quad \mbox{\rm weakly in } L^2(0,T;H^{-1}(\mathbb{R}^d)) .$$
Let us prove that actually
$$\rho ^{0,v_n}\longrightarrow \rho ^{0,v^*} \quad \mbox{\rm strongly in } L^2(Q_T) .$$
The first step is to show that
$$\lim_{\varepsilon \rightarrow 0}\int_0^T\int_{B(0_d;\sqrt{\frac{2}{\varepsilon }})^c}\rho ^{0,v_n}(t,x)^2dx \ dt =0 \quad \mbox{\rm uniformly in } n \, . $$
Indeed, let $\psi \in C^2([0,+\infty ])$ given by
$$\psi (r)=\left\{ \begin{array}{ll}
	0, & r\in [0,1], \\
	6(r-1)^5-15(r-1)^4+10(r-1)^3, \quad & r\in (1,2), \\
	1, & r\in [2, +\infty ) 
\end{array}
	\right. $$  
and $\psi_{\varepsilon}(x)=\psi (\varepsilon |x|^2) $. 
It follows that
$$\frac{\partial \psi _{\varepsilon }}{\partial x_i}(x)=2\varepsilon x_i\psi '(\varepsilon |x|^2), $$
$$\frac{\partial ^2\psi _{\varepsilon }}{\partial x_i\partial x_j}(x)=4\varepsilon ^2x_ix_j\psi ''(\varepsilon |x|^2)
+2\varepsilon \delta _{ij}\psi '(\varepsilon |x|^2) .$$
We have that $\psi _{\varepsilon}\rho ^{0,v_n}$ satisfies
\begin{equation*}
	\left\{ \begin{array}{lll}
		\displaystyle
		\frac{\partial (\psi_{\varepsilon }\rho^{0,v_n})}{\partial t}& \displaystyle
		=-\nabla \cdot (v_n\psi _{\varepsilon}\rho^{0,v_n})+\frac{\partial ^2(a_{ij}\psi_{\varepsilon }\rho ^{0,v_n})}{\partial x_i\partial x_j}, & \vspace{2mm} \\
		~ & \ \ \ \displaystyle +v_n\rho ^{0,v_n}\cdot \nabla \psi _{\varepsilon}-2\frac{\partial (a_{ij}\rho ^{0,v_n})}{\partial x_i}\frac{\partial \psi _{\varepsilon }}{\partial x_j}-a_{ij}\rho ^{0,v_n}\frac{\partial ^2\psi _{\varepsilon }}{\partial x_i\partial x_j}, \quad & (t,x)\in Q_T, \vspace{2mm} \\
		\psi_{\varepsilon}(x)\rho ^{0,v_n}(0,x)& =\psi _{\varepsilon }(x)\rho _0(x), & x\in \mathbb{R}^d .
	\end{array}
	\right.
\end{equation*}
It follows that
$$\frac{1}{2}\| \psi _{\varepsilon }\rho ^{0,v_n}(t)\|^2_{L^2}-\frac{1}{2}\| \psi _{\varepsilon }\rho _0\|^2_{L^2}
=\int_0^t\int_{\mathbb{R}^d}v_n(s,x)\psi_{\varepsilon }(x)\rho ^{0,v_n}(s,x)\cdot \nabla (\psi_{\varepsilon}(x)\rho ^{0,v_n}(s,x))dx \ ds $$
$$-\int_0^t\int_{\mathbb{R}^d}\Big[ a_{ij}\frac{\partial (\psi_{\varepsilon }\rho ^{0,v_n})}{\partial x_j}\frac{\partial (\psi_{\varepsilon }\rho ^{0,v_n})}{\partial x_i}
+\frac{\partial a_{ij}}{\partial x_j}\psi_{\varepsilon }\rho ^{0,v_n}\frac{\partial (\psi_{\varepsilon }\rho ^{0,v_n}) }{\partial x_i}\Big] dx \ ds $$
$$+\int_0^t\int_{\mathbb{R}^d}\Big[ v_n\rho ^{0,v_n}\cdot \nabla \psi_{\varepsilon}(\psi_{\varepsilon }\rho^{0,v_n})
+2a_{ij}\rho ^{0,v_n}\frac{\partial }{\partial x_i}(\frac{\partial \psi_{\varepsilon }}{\partial x_j}\psi_{\varepsilon }\rho ^{0,v_n})
-a_{ij}\rho ^{0,v_n}\frac{\partial ^2\psi_{\varepsilon }}{\partial x_i\partial x_j}\psi_{\varepsilon }\rho^{0,v_n} \Big] dx \ ds . $$
After an easy calculation and using the properties of $\sigma $ and $\psi $ we get that there exists a positive constant $C_0$ such that
$$\| \psi _{\varepsilon }\rho ^{0,v_n}(t)\|^2_{L^2}+\int_0^t\int_{\mathbb{R}^d}|\nabla (\psi_{\varepsilon }\rho^{0,v_n})|^2dx \ ds\leq C_0(\| \psi _{\varepsilon }\rho _0\|^2_{L^2}+\varepsilon) , \quad \forall \varepsilon >0, \forall t\in [0,T] .$$
It follows that 
$$\lim_{\varepsilon \rightarrow 0}\int_0^T\int_{B(0_d;\sqrt{\frac{2}{\varepsilon }})^c}\rho ^{0,v_n}(t,x)^2dx \ dt =0 \quad \mbox{\rm uniformly in } n \, .$$
We may conclude that
$$\int_{Q_T}|\rho^{0,v_n}(t,x)-\rho ^{0,v^*}(t,x)|^2dx \ dt=\int_0^T\int_{B(0_d;\sqrt{\frac{2}{\varepsilon }})}|\rho ^{0,v_n}(t,x)-\rho ^{0,v^*}|^2dx \ dt$$
$$+\int_0^T\int_{B(0_d;\sqrt{\frac{2}{\varepsilon }})^c}|\rho ^{0,v_n}(t,x)-\rho ^{0,v^*}(t,x)|^2dx \ dt\leq \int_0^T\int_{B(0_d;\sqrt{\frac{2}{\varepsilon }})}|\rho ^{0,v_n}(t,x)-\rho ^{0,v^*}|^2dx \ dt$$
$$+\int_0^T\int_{B(0_d;\sqrt{\frac{2}{\varepsilon }})^c}\rho ^{0,v_n}(t,x)^2dx \ dt+\int_0^T\int_{B(0_d;\sqrt{\frac{2}{\varepsilon }})^c}\rho ^{0,v^*}(t,x)^2dx \ dt $$
and consequently we get that
$$\rho ^{0,v_n}\longrightarrow \rho ^{0,v^*} \quad \mbox{\rm strongly in } L^2(Q_T).$$
We conclude (using assumptions (H6) or (H6)') that there exists $u^*\in {\cal U}$ such that
$$v^*(t,x)=f(t,x,u^*(t,x)), \quad \mbox{\rm a.e. } (t,x)\in Q_T$$
and consequently 
$$\rho ^{u_n}\longrightarrow \rho ^{u^*} \quad \mbox{\rm strongly in } L^2(Q_T).$$
\vspace{2mm}

The properties of $G$ and $G^0$ imply that
$$\liminf _{n\rightarrow +\infty }\Big\{ \int_{Q_T}G(t,x,\rho ^{u_n}(t,x))dx \ dt+\int_{\mathbb{R}^d}G^0(x,\rho ^{u_n}(T,x))dx\Big\} $$
$$\geq \int_{Q_T} G(t,x,\rho ^{u^*}(t,x))dx \ dt+\int_{\mathbb{R}^d}G^0(x,\rho ^{u^*}(T,x))dx \, .$$
On the other hand if $f$ satisfies (H6) we may choose $u^*\in {\cal U}$ (satisfying $v^*(t,x)=f(t,x,u^*(t,x))$ a.e. $(t,x)\in Q_T$) such that $u_n\longrightarrow u^*$ weak * in $L^{\infty }(Q_T)^{\ell}$ and so
$$\liminf _{n\rightarrow +\infty }\int_{Q_T}L(t,x,u_n(t,x))\rho^{u_n}(t,x)dx \ dt 
=\liminf _{n\rightarrow +\infty }\int_{Q_T}L(t,x,u_n(t,x))\rho^{u^*}(t,x)dx \ dt $$
$$\geq \int_{Q_T}L(t,x,u^*(t,x))\rho^{u^*}(t,x)dx \ dt .$$
We may infer that $u^*$ is an optimal control for the optimal control problem.
\vspace{2mm}

If $f$ satisfies (H6)', then 
$$\liminf _{n\rightarrow +\infty }\int_{Q_T}L(t,x,u_n(t,x))\rho^{u_n}(t,x)dx \ dt= \liminf _{n\rightarrow +\infty }\int_{Q_T}L(t,x,f^{-1}(t,x,v_n(t,x)))\rho^{u_n}(t,x)dx \ dt$$
$$=\liminf _{n\rightarrow +\infty }\int_{Q_T}L(t,x,f^{-1}(t,x,v_n(t,x)))\rho^{u^*}(t,x)dx \ dt $$
$$\geq \int_{Q_T}L(t,x,f^{-1}(t,x,v^*(t,x)))\rho^{u^*}(t,x)dx \ dt=\int_{Q_T}L(t,x,u^*(t,x))\rho^{u^*}(t,x)dx \ dt . $$
We may conclude that $u^*$ is an optimal control for the optimal control problem.
\vspace{5mm}

Consider the following approximating optimal control problem

$$\underset{u \in {\cal U}}{\text{Minimize}} \ \Big\{ J_{\varepsilon }(u,\rho ^u)+\frac{1}{2}\int_{Q_T}\varphi (x)\rho ^u(t,x)|u(t,x)-u^*(t,x)|^2 dx \, dt\Big\} ,$$
where $\varphi \in C_b(\mathbb{R}^d)\cap L^1(\mathbb{R}^d)$, $\varphi (x)>0$, $\forall x\in \mathbb{R}^d$, $\varepsilon >0$ and 
$$J_{\varepsilon }(u,\rho ^u)=\int_{Q_T}[L(t,x,u(t,x))\rho ^u(t,x)+G_{\varepsilon }(t,x,\rho ^u(t,x))] dx \, dt+\int_{\mathbb{R}^d}
G^0_{\varepsilon }(x,\rho ^u(T,x))dx \, .$$
Here
$$G_{\varepsilon }(t,x,r)=\min \Big\{ \frac{1}{2\varepsilon }|r-\theta |^2 +G(t,x,\theta ); \ \theta \in \mathbb{R}\Big\} , \quad r\in \mathbb{R}, $$
$$G^0_{\varepsilon }(x,r)=\min \Big\{ \frac{1}{2\varepsilon }|r-\theta |^2 +G^0(x,\theta ); \ \theta \in \mathbb{R}\Big\} , \quad r\in \mathbb{R}, $$
are the Yosida approximations of $r\mapsto G(t,x,r)$ and $r\mapsto G^0(x,r)$, respectively. It is well known that
$G_{\varepsilon }$ and $G^0_{\varepsilon }$ are differentiable on $\mathbb{R}$ a.e. $(t,x)\in Q_T$, $x\in \mathbb{R}^d$, respectively
and that
$$G_{\varepsilon }(t,x,r)=\frac{1}{2\varepsilon }|r-(I+\varepsilon G_r(t,x,\cdot ))^{-1}r|^2+G(t,x,(I+\varepsilon G_r(t,x,\cdot ))^{-1}r),$$
$$G^0_{\varepsilon }(x,r)=\frac{1}{2\varepsilon }|r-(I+\varepsilon G^0_r(x,\cdot ))^{-1}r|^2+G^0(x,(I+\varepsilon G^0_r(x,\cdot ))^{-1}r),$$
for any $r\in \mathbb{R}$. We shall denote the derivatives with respect to $r$ by $G_{\varepsilon}'(t,x,r)$ and $(G^0_{\varepsilon })'(x,r)$, respectively.

Let $(u_{\varepsilon },\rho _{\varepsilon })$ be an optimal pair for the above mentioned optimal control problem. It is clear that
$\rho _{\varepsilon }$ satisfies (is a weak solution to)
\begin{equation}\label{eq1bis2}
	\left\{ \begin{array}{ll}
		\displaystyle
		\frac{\partial \rho_{\varepsilon }}{\partial t}(t,x)=-\nabla \cdot (f^{u_{\varepsilon }}(t,x)\rho_{\varepsilon }(t,x))+\frac{\partial ^2(a_{ij}\rho_{\varepsilon })}{\partial x_i\partial x_j}(t,x), \quad &(t,x)\in Q_T,\vspace{2mm} \\
		\rho_{\varepsilon }(0,x)=\rho_0(x), & x\in \mathbb{R}^d.
	\end{array}
	\right.
\end{equation}
There exists a unique weak solution to (\ref{eq1bis2}) and we also get that
$$\frac{1}{2}\| \rho _{\varepsilon }(t)\|^2_{L^2}=\int_0^t\int_{\mathbb{R}^d}f^{u_{\varepsilon }}(s,x)\rho _{\varepsilon }(s,x)\cdot \nabla \rho _{\varepsilon }(s,x)dx \, ds$$
$$-\int_0^t\int_{\mathbb{R}^d}\frac{\partial  (a_{ij}\rho_{\varepsilon })}{\partial x_j}(s,x)\frac{\partial \rho _{\varepsilon }}{\partial x_i}(s,x)dx \, ds
+\frac{1}{2}\|\rho _0\|^2_{L^2} .$$
It follows after an easy evaluation that there exists a positive constant $C_{01}$ such that
$$\| \rho _{\varepsilon }(t)\|^2_{L^2}+\int_0^t\int_{\mathbb{R}^d}| \nabla \rho _{\varepsilon }(s,x)|^2dx \, ds+\Big\| \frac{d\rho _{\varepsilon }}{dt}\Big\|^2_{L^2(0,T;H^{-1})}\leq C_{01}\| \rho _0\|^2_{L^2},$$
for any $t\in [0,T]$, $\forall \varepsilon >0$. 

Let $v\in L^{\infty }(Q_T)^{\ell }\cap L^1(Q_T)^{\ell }$. For any such $v$ and $\lambda >0$ we have that
$$J_{\varepsilon }(u_{\varepsilon }+\lambda v,\rho ^{u_{\varepsilon }+\lambda v})+\frac{1}{2}\int_{Q_T}\varphi \rho ^{u_{\varepsilon }+\lambda v}
|u_{\varepsilon }+\lambda v-u^*|^2dx \, dt
\geq J_{\varepsilon }(u_{\varepsilon },\rho ^{u_{\varepsilon }})+\frac{1}{2}\int_{Q_T}\varphi \rho ^{u_{\varepsilon }}|u_{\varepsilon }-u^*|^2dx \, dt.$$
It follows that
$$\int_{Q_T}\Big[ L'(t,x,u_{\varepsilon }(t,x),v(t,x))\rho^{u_{\varepsilon }}(t,x)+L(t,x,u_{\varepsilon }(t,x))z(t,x)+G_{\varepsilon }'(t,x,\rho _{\varepsilon }(t,x))z(t,x)$$
$$+\varphi (x)\rho _{\varepsilon }(t,x)(u_{\varepsilon }(t,x)-u^*(t,x))\cdot v(t,x)+\frac{1}{2}\varphi (x)|u_{\varepsilon }(t,x)-u^*(t,x)|^2z(t,x)\Big]
dx \, dt$$
\begin{equation}\label{ineq}
+\int_{\mathbb{R}^d}(G^0_{\varepsilon })'(x,\rho _{\varepsilon }(T,x))z(T,x)dx \geq 0, 
\end{equation}
where $\displaystyle L'(t,x,u,v)=\lim_{\theta \rightarrow 0+}\frac{L(t,x,u+\theta v)-L(t,x,u)}{\theta }$ and $z$ is the weak solution to
\begin{equation}\label{eqzbis}
	\left\{ \begin{array}{ll}
		\displaystyle
		\frac{\partial z}{\partial t}(t,x)=-\nabla \cdot (f^{u_{\varepsilon }}(t,x)z(t,x)+f_u^{u_{\varepsilon }}(t,x)v(t,x)\rho _{\varepsilon }(t,x)) & \vspace{2mm} \\
\hspace{3cm} \displaystyle+\frac{\partial ^2(a_{ij}z)}{\partial x_i\partial x_j}(t,x), \quad &(t,x)\in Q_T, \vspace{2mm} \\
		z(0,x)=0, & x\in \mathbb{R}^d.
	\end{array}
	\right.
\end{equation}
Let $p_{\varepsilon }$ be the weak solution to
\begin{equation}\label{dual}
	\left\{ \begin{array}{ll}
		\displaystyle
		\frac{\partial p_{\varepsilon }}{\partial t}(t,x)=-f^{u_{\varepsilon }}(t,x)\cdot \nabla p_{\varepsilon }(t,x)
		-a_{ij}(t,x)\frac{\partial ^2p_{\varepsilon }}{\partial x_i\partial x_j}(t,x) \vspace{2mm} \\
		\hspace{2cm} +G_{\varepsilon }'(t,x,\rho^{u_{\varepsilon }}(t,x))+L(t,x,u_{\varepsilon }(t,x)) \vspace{2mm} \\
		\hspace{2cm} \displaystyle +\frac{1}{2}\varphi (x)|u_{\varepsilon }(t,x)-u^*(t,x)|^2, \quad &(t,x)\in Q_T,\vspace{2mm} \\
		p_{\varepsilon }(T,x)=-(G^0_{\varepsilon })'(x,\rho _{\varepsilon }(T,x)), & x\in \mathbb{R}^d.
	\end{array}
	\right.
\end{equation}

\noindent
Let us establish some evaluations concerning $g_{\varepsilon }, (g_0)_{\varepsilon }$ and their derivatives.

We get that
$$G_{\varepsilon }(t,x,r)\leq \min \Big\{ \frac{1}{2\varepsilon }|\theta -r|^2+\bar{C}_1\theta ^2+\bar{\alpha }(t,x); \ \theta\in \mathbb{R}\Big\} $$
(and since the the minimum is attained for $\theta _{min}=\frac{r}{1+2\bar{C}_1\varepsilon}$)
$$= \frac{1}{2\varepsilon }\Big| \frac{r}{1+2\bar{C}_1\varepsilon }-r\Big|^2+\bar{C}_1\frac{r^2}{(1+2\bar{C}_1\varepsilon)^2}+\bar{\alpha }(t,x)$$
(which gives after an easy calculation)
$$=\frac{\bar{C}_1}{1+2\bar{C}_1\varepsilon }r^2+\bar{\alpha }(t,x) $$
for any $(t,x,r)\in Q_T\times \mathbb{R}$.

On the other hand we get that
$$G_{\varepsilon }(t,x,r)\geq \min \Big\{ \frac{1}{2\varepsilon }|\theta -r|^2+\alpha(t,x) \theta ; \ \theta\in \mathbb{R}\Big\} $$
(and since the the minimum is attained for $\theta _{min}=r-\varepsilon \alpha (t,x)$)
$$= \frac{1}{2\varepsilon }\varepsilon ^2\alpha (t,x)^2+\alpha (t,x)(r-\varepsilon \alpha (t,x))$$
$$=r \, \alpha (t,x)-\frac{\varepsilon }{2}\alpha (t,x)^2$$
for any $(t,x,r)\in Q_T\times \mathbb{R}$.

Since for any $(t,x)\in Q_T$, the mapping $r\mapsto G_{\varepsilon }(t,x,r)$ is convex we may infer that
for any $(t,x,r)\in Q_T\times \mathbb{R}$:
$$G_{\varepsilon }'(t,x,r)h\leq G_{\varepsilon }(t,x,r+h)-G_{\varepsilon }(t,x,r) $$
$$\leq \frac{\bar{C}_1}{1+2\bar{C}_1\varepsilon }(r+h)^2+\bar{\alpha }(t,x)-r \, \alpha (t,x)+\frac{\varepsilon }{2}\alpha (t,x)^2$$
$$=\frac{\bar{C}_1}{1+2\bar{C}_1\varepsilon }h^2 +\frac{2\bar{C}_1r}{1+2\bar{C}_1\varepsilon }h+\frac{\bar{C}_1}{1+2\bar{C}_1\varepsilon }r^2
+\bar{\alpha }(t,x)-r \, \alpha (t,x)+\frac{\varepsilon }{2}\alpha (t,x)^2.$$
It follows that for any $(t,x,r)\in Q_T\times \mathbb{R}$:
$$\frac{\bar{C}_1}{1+2\bar{C}_1\varepsilon }h^2 +\Big(
\frac{2\bar{C}_1r}{1+2\bar{C}_1\varepsilon }-G_{\varepsilon }'(t,x,r)\Big)h+\frac{\bar{C}_1}{1+2\bar{C}_1\varepsilon }r^2
+\bar{\alpha }(t,x)-r \, \alpha (t,x)+\frac{\varepsilon }{2}\alpha (t,x)^2\geq 
0, \ \forall h\in \mathbb{R} .$$
This happens if and only if
$$\Big|\frac{2\bar{C}_1r}{1+2\bar{C}_1\varepsilon }-G_{\varepsilon }'(t,x,r)\Big|^2-
\frac{4\bar{C}_1}{1+2\bar{C}_1\varepsilon }
\Big( \frac{\bar{C}_1}{1+2\bar{C}_1\varepsilon }r^2
+\bar{\alpha }(t,x)-r \, \alpha (t,x)+\frac{\varepsilon }{2}\alpha (t,x)^2\Big)\leq 0 .$$
It follows that
$$\Big|\frac{2\bar{C}_1r}{1+2\bar{C}_1\varepsilon }-G_{\varepsilon }'(t,x,r)\Big| \leq 
\sqrt{4\bar{C}_1	\Big( \bar{C}_1r^2	+\bar{\alpha }(t,x)+\frac{r^2}{2}+\frac{1}{2}\alpha (t,x)^2+\frac{\varepsilon }{2}\alpha (t,x)^2\Big)}
$$
for any $(t,x,r)\in Q_T\times \mathbb{R}$. We may conclude that for any $\varepsilon \in (0,\varepsilon _0)$:
$$|G_{\varepsilon }'(t,x,r)|\leq M_G|r|+A(t,x), \quad \forall (t,x,r)\in Q_T\times \mathbb{R},$$
where $M_G$ is a positive constant and $A\in L^1(Q_T)\cap L^{\infty }(Q_T)$.
\vspace{5mm}

Analogously we get that there exists a positive constant $M_{G^0}$ and a function $A^0\in L^1(\mathbb{R}^d)\cap L^{\infty }(\mathbb{R}^d)$
such that
$$|(G^0_{\varepsilon })'(x,r)|\leq M_{G^0}|r|+A^0(x), \quad \forall (x,r)\in \mathbb{R}^d\times \mathbb{R}.$$
Multiplying (\ref{dual}) by $z$ and integrating over $Q_T$ we obtain that
$$-\int_{\mathbb{R}^d}(G^0_{\varepsilon})'(x,\rho_{\varepsilon}(T,x))z(T,x)dx-\int_{Q_T}\Big[\nabla p_{\varepsilon }(t,x)\cdot
(f^{u_{\varepsilon }}(t,x)u_{\varepsilon}(t,x)z(t,x) $$
$$+f_u^{u_{\varepsilon }}(t,x)v(t,x)\rho _{\varepsilon}(t,x))-\frac{\partial p_{\varepsilon }}{\partial x_i}(t,x)\frac{\partial (a_{ij}z)}{\partial x_j}(t,x)\Big] dx \, dt$$
$$=\int_{Q_T}\Big[ -z(t,x)f^{u_{\varepsilon }}(t,x)\cdot \nabla p_{\varepsilon }(t,x)+\frac{\partial (a_{ij}z)}{\partial x_i}(t,x)\frac{\partial p_{\varepsilon }}{\partial x_j}(t,x)$$
$$+G'_{\varepsilon }(t,x,\rho _{\varepsilon }(t,x))z(t,x)+L(t,x,u_{\varepsilon }(t,x))z(t,x)+\frac{1}{2}\varphi (x)|u_{\varepsilon }(t,x)-u^*(t,x)|^2z(t,x)\Big] dx \, dt.$$
It follows that
$$\int_{Q_T}\Big[ L(t,x,u_{\varepsilon }(t,x))z(t,x)+G_{\varepsilon }'(t,x,\rho _{\varepsilon }(t,x))z(t,x)+\frac{1}{2}\varphi (x)|u_{\varepsilon }(t,x)-u^*(t,x)|^2z(t,x)\Big] dx \, dt$$
$$+\int_{\mathbb{R}^d}(G^0_{\varepsilon })'(x,\rho_{\varepsilon }(T,x))z(T,x)dx =-\int_{Q_T}\nabla p_{\varepsilon }(t,x)\cdot f_u^{u_{\varepsilon }}(t,x)v(t,x)\rho _{\varepsilon }(t,x)dx \, dt .$$
By (\ref{ineq}) we obtain that
$$\int_{Q_T}\rho _{\varepsilon }(t,x)\Big[ L'(t,x,u_{\varepsilon }(t,x),v(t,x))+\varphi (x)(u_{\varepsilon }(t,x)-u^*(t,x))\cdot v(t,x)$$
$$-f_u^{u_{\varepsilon }}(t,x)^T \nabla p_{\varepsilon }(t,x)\cdot v(t,x) \Big]dx \, dt \geq 0$$
and so we may conclude that
$$f_u^{u_{\varepsilon }}(t,x)^T \nabla p_{\varepsilon }(t,x)-\varphi (x)(u_{\varepsilon }(t,x)-u^*(t,x))\in L_u(t,x,u_{\varepsilon }(t,x))$$
for any $(t,x)\in Q_T$.

It follows that for almost any $(t,x)\in Q_T$
$$L(t,x,u_{\varepsilon }(t,x))+H(t,x,f_u^{u_{\varepsilon }}(t,x)^T \nabla p_{\varepsilon }(t,x)-\varphi (x)(u_{\varepsilon }(t,x)-u^*(t,x)))$$
$$=u_{\varepsilon }(t,x)\cdot (f_u^{u_{\varepsilon }}(t,x)^T \nabla p_{\varepsilon }(t,x)-\varphi (x)(u_{\varepsilon }(t,x)-u^*(t,x))) $$
and so $p_{\varepsilon }$ is the weak solution to 
\begin{equation*}
	\left\{ \begin{array}{ll}
\displaystyle
\frac{\partial p_{\varepsilon }}{\partial t}(t,x)& \displaystyle +a_{ij}(t,x)\frac{\partial ^2p_{\varepsilon }}{\partial x_i\partial x_j}(t,x) \vspace{2mm} \\
~ & \ \ +H(t,x,f_u^{u_{\varepsilon }}(t,x)^T \nabla p_{\varepsilon }(t,x)-\varphi (x)(u_{\varepsilon }(t,x)-u^*(t,x))) \vspace{2mm} \\
~ & = u_{\varepsilon }(t,x)\cdot (f_u^{u_{\varepsilon }}(t,x)^T \nabla p_{\varepsilon }(t,x)-\varphi (x)(u_{\varepsilon }(t,x)-u^*(t,x))) \vspace{2mm} \\
~ &  \ \ \displaystyle -f^{u_{\varepsilon }}(t,x)\cdot \nabla p_{\varepsilon }(t,x)
+G_{\varepsilon }'(t,x,\rho_{\varepsilon }(t,x))+\frac{1}{2}\varphi (x)|u_{\varepsilon }(t,x)-u^*(t,x)|^2 \vspace{2mm} \\
~ & = (f_u^{u_{\varepsilon }}(t,x)u_{\varepsilon }(t,x)-f^{u_{\varepsilon }}(t,x))\cdot \nabla p_{\varepsilon }(t,x) \vspace{2mm} \\
~ &  \ \ \displaystyle +G_{\varepsilon }'(t,x,\rho_{\varepsilon }(t,x))-\frac{1}{2}\varphi (x)(|u_{\varepsilon }(t,x)|^2-|u^*(t,x)|^2),
\quad (t,x)\in Q_T\vspace{2mm} \\
		p_{\varepsilon }(T,x)&=-(G^0_{\varepsilon })'(x,\rho _{\varepsilon }(T,x)), \quad x\in \mathbb{R}^d.
	\end{array}
	\right.
\end{equation*}
On the other hand 
$$u_{\varepsilon }(t,x)=H_r(t,x,f_u^{u_{\varepsilon }}(t,x)^T \nabla p_{\varepsilon }(t,x)-\varphi (x)(u_{\varepsilon }(t,x)-u^*(t,x))), \quad \mbox{\rm a.e. } (t,x)\in Q_T .$$
We have that
$$J_{\varepsilon }(u_{\varepsilon },\rho _{\varepsilon})+\frac{1}{2}\int_{Q_T}\varphi (x)\rho _{\varepsilon}(t,x)|u_{\varepsilon }(t,x)-u^*(t,x)|^2dx \, dt \leq J_{\varepsilon }(u^*,\rho ^*) $$
and
$$J_{\varepsilon }(u^*,\rho ^*)=\int_{Q_T}[L(t,x,u^*(t,x))\rho ^*(t,x)+G_{\varepsilon }(t,x,\rho ^*(t,x))] dx \, dt
+\int_{\mathbb{R}^d}G^0_{\varepsilon }(x,\rho ^*(T,x))dx \, .$$
Since 
$$G(t,x,\rho^*(t,x))\geq G_{\varepsilon }(t,x,\rho^*(t,x))\longrightarrow  G(t,x,\rho^*(t,x)),\quad \mbox{\rm a.e. } (t,x)\in Q_T,$$
$$G^0(x,\rho^*(T,x))\geq G^0_{\varepsilon }(x,\rho^*(T,x))\longrightarrow  G^0(x,\rho^*(T,x)),\quad \mbox{\rm a.e. } x\in \mathbb{R}^d$$
we conclude that
$$\lim_{\varepsilon \rightarrow 0}J_{\varepsilon }(u^*,\rho ^*)=J(u^*,\rho ^*) \, .$$
\vspace{3mm}

We have that on a subsequence (also indexed by $\varepsilon \rightarrow 0$):
$$f^{u_{\varepsilon }}\longrightarrow \tilde{v} \quad \mbox{\rm weak * in } L^{\infty }(Q_T)^d , $$
$$\rho _{\varepsilon }\longrightarrow \rho ^{0,\tilde{v}} \quad \mbox{\rm weakly in } L^2(0,T; H^1(\mathbb{R}^d)), $$
and via Aubin's compactness theorem, also strongly in $L^2(0,T; L^2_{loc}(\mathbb{R}^d))$, 
$$\rho _{\varepsilon }(T)\longrightarrow \rho ^{0,\tilde{v}}(T) \quad \mbox{\rm weakly in } L^2(\mathbb{R}^d) .$$
Actually, it is possible to prove (as in the first part of the proof) that
$$\rho _{\varepsilon }\longrightarrow \rho ^{0,\tilde{v}} \quad \mbox{\rm strongly in } L^2(0,T; L^2(\mathbb{R}^d)) \, . $$
Since $f$ satisfies (H6) or (H6)' we get that there exists $\tilde{u}\in {\cal U}$ such that $f(t,x,\tilde{u}(t,x))=\tilde{v}(t,x)$ a.e. $(t,x)\in Q_T$ and so
$$\rho _{\varepsilon }\longrightarrow \rho ^{\tilde{u}} \quad \mbox{\rm strongly in } L^2(0,T; L^2(\mathbb{R}^d)) \, . $$
It follows using a previous argument that
$$\liminf_{\varepsilon \rightarrow 0}\Big\{ \int_{Q_T}[L(t,x,u_{\varepsilon }(t,x))\rho _{\varepsilon }(t,x)+G_{\varepsilon }(t,x,\rho _{\varepsilon }(t,x))] dx \, dt\Big\} $$
$$ \geq \int_{Q_T}[L(t,x,\tilde{u}(t,x))\rho ^{\tilde{u}}(t,x)+G(t,x,\rho ^{\tilde{u}}(t,x))] dx \, dt .$$

On the other hand we have that for an arbitrary but fixed $\varepsilon _0>0$ and for any $0<\varepsilon <\varepsilon _0$:
$$\int_{\mathbb{R}^d}G^0_{\varepsilon }(x,\rho _{\varepsilon }(T,x))dx\geq \int_{\mathbb{R}^d}G^0_{\varepsilon _0}(x,\rho _{\varepsilon }(T,x))dx,$$
and so
$$\liminf_{\varepsilon \rightarrow 0}\int_{\mathbb{R}^d}G^0_{\varepsilon }(x,\rho _{\varepsilon }(T,x))dx\geq \int_{\mathbb{R}^d}G^0_{\varepsilon _0}(x,\rho ^{\tilde{u}}(T,x))dx,$$
for any $\varepsilon _0>0$. If we apply again Fatou's lemma we obtain that
$$\liminf_{\varepsilon \rightarrow 0}\int_{\mathbb{R}^d}G^0_{\varepsilon }(x,\rho _{\varepsilon }(T,x))dx\geq \int_{\mathbb{R}^d}G^0(x,\rho ^{\tilde{u} }(T,x))dx$$
and so
$$\liminf_{\varepsilon \rightarrow 0}\Big\{ \int_{Q_T}[L(t,x,u_{\varepsilon }(t,x))\rho _{\varepsilon }(t,x)+G_{\varepsilon }(t,x,\rho _{\varepsilon }(t,x))] dx \, dt
+\int_{\mathbb{R}^d}G^0_{\varepsilon }(x,\rho _{\varepsilon }(T,x))dx\Big\} $$
$$ \geq \int_{Q_T}[L(t,x,\tilde{u}(t,x))\rho ^{\tilde{u}}(t,x)+G(t,x,\rho ^{\tilde{u}}(t,x))] dx \, dt
+\int_{\mathbb{R}^d}G^0(x,\rho ^{\tilde{u}}(T,x))dx \, .$$

We may conclude that
$$J(u^*,\rho ^*)\leq J(\tilde{u},\rho ^{\tilde{u}})+\limsup_{\varepsilon \rightarrow 0}\frac{1}{2}\int_{Q_T}\varphi (x)\rho_{\varepsilon }(t,x)|u_{\varepsilon}(t,x)-u^*(t,x)|^2dx \, dt
\leq J(u^*,\rho ^*) \, $$
and it follows that
$$\lim_{\varepsilon \rightarrow 0}\int_{Q_T}\varphi (x)\rho_{\varepsilon }(t,x)|u_{\varepsilon}(t,x)-u^*(t,x)|^2dx \, dt = 0.$$
We may conclude that on a subsequence
$$u_{\varepsilon }(t,x)\longrightarrow u^*(t,x), \quad \mbox{\rm a.e. } (t,x)\in Q_T \, .$$
\vspace{3mm}

\noindent
\bf Lemma 3.2. \it If $\varepsilon _n \rightarrow 0$, then there exists a subsequence (also denoted by $\{ \varepsilon_n\}$)
such that
$$G_{\varepsilon _n}'(t,x,\rho _{\varepsilon _n}(t,x))\longrightarrow \eta (t,x), \quad \mbox{\rm weakly in } L^2(Q_T),$$
$$(G^0_{\varepsilon _n})'(x,\rho _{\varepsilon _n}(T,x))\longrightarrow \eta _0(x), \quad \mbox{\rm weakly in } L^2(\mathbb{R}^d) \, .$$

Proof. \rm Indeed,
$$|G_{\varepsilon _n}'(t,x,\rho _{\varepsilon _n}(t,x))|\leq M_G\rho _{\varepsilon _n}(t,x)+A(t,x)$$
and since $\{ M_G\rho_{\varepsilon _n}+A\}$ is bounded in $L^2(Q_T)$ we get the first conclusion. 

The second convergence follows analogously.
\vspace{3mm}

Multiplying (\ref{dual}) by $p_{\varepsilon _n}$ and since $\{ G_{\varepsilon _n}'(t,x,\rho _{\varepsilon _n}(t,x))\}_{n\in \mathbb{N}^*} $, 
$\{ L(t,x, u_{\varepsilon _n}(t,x))\}_{n\in \mathbb{N}^*} $, \break
$\{ \frac{1}{2}\varphi (x)(u_{\varepsilon _n}(t,x)-u^*(t,x))\}_{n\in \mathbb{N}^*} $ are bounded in $L^2(Q_T)$ and 
$\{ (G^0_{\varepsilon_n})'(x,\rho _{\varepsilon _n}(T,x))\}_{n\in \mathbb{N}^*} $ is bounded in $L^2(\mathbb{R}^d)$ it follows that $\{ p_{\varepsilon _n}\}_{n\in \mathbb{N}^*}$
and $\{ \nabla p_{\varepsilon _n}\}_{n\in \mathbb{N}^*}$ are bounded in $L^2(Q_T)$ and $L^2(Q_T)^d$, respectively.

It follows that on a subsequence 
$$p_{\varepsilon _n}\longrightarrow p \quad \mbox{\rm weakly in } L^2(Q_T),$$
$$\nabla p_{\varepsilon _n}\longrightarrow \nabla p \quad \mbox{\rm weakly in } L^2(Q_T)^d,$$
$$H(\cdot ,\cdot ,(f_u^{u_{\varepsilon _n}})^T\nabla p_{\varepsilon _n}-\varphi (u_{\varepsilon _n}-u^*))\longrightarrow \zeta ^* \quad \mbox{\rm weakly in } L^2(Q_T) \, .$$

On the other hand 
$$\rho _{\varepsilon _n}\longrightarrow \rho ^* \quad \mbox{\rm in } C([0,T]; L^2(\mathbb{R}^d)) \, .$$
Since for any $r\in \mathbb{R}$:
$$(G^0_{\varepsilon _n})'(x,\rho _{\varepsilon _n}(T,x))(r-\rho _{\varepsilon _n}(T,x))\leq G^0_{\varepsilon _n}(x,r)-G^0_{\varepsilon _n}(x,\rho _{\varepsilon _n}(T,x))$$
it follows that
$$\eta_0(x)(r-\rho _{\varepsilon }(T,x))\leq G^0(x,r)-G^0(x,\rho ^*(T,x))$$
and so 
$$\eta _0(x)\in G^0_r(x,\rho ^*(T,x)), \quad \mbox{\rm a.e. } x\in \mathbb{R}^d \, .$$
In the same manner it follows that
$$\eta (t,x)\in G_r(t,x,\rho ^*(t,x)), \quad \mbox{\rm a.e. } (t,x)\in Q_T \, .$$

Since
$$u_{\varepsilon _n}(t,x)\in H_u(t,x,f_u^{u_{\varepsilon _n}}(t,x)^T\nabla p_{\varepsilon _n}(t,x)-\varphi (x)(u_{\varepsilon _n}(t,x)-u^*(t,x)))$$
it follows that for any $\theta \in \mathbb{R}^{\ell }$:
$$\int_{Q_T}u_{\varepsilon _n}(t,x)\cdot (\theta -f_u^{u_{\varepsilon _n}}(t,x)^T\nabla p_{\varepsilon _n}(t,x)-\varphi (x)(u_{\varepsilon _n}-u^*)(t,x))\psi (t,x)dx \, dt$$
$$\leq \int_{Q_T}(H(t,x,\theta )-H(t,x,f_u^{u_{\varepsilon _n}}(t,x)^T\nabla p_{\varepsilon _n}(t,x)-\varphi (x)(u_{\varepsilon _n}-u^*)))\psi (t,x)dx \, dt$$
for any $\psi \in L^{\infty }(Q_T)$ with compact support and $\psi (t,x)\geq 0$.
We get that
$$\int_{Q_T}u^*(t,x)\cdot (\theta -f_u^{u^*}(t,x)^T\nabla p(t,x))\psi (t,x)dx \, dt$$
$$\leq \int_{Q_T}(H(t,x,\theta )-H(t,x,f_u^{u^*}(t,x)^T\nabla p(t,x))\psi (t,x)dx \, dt$$
for any $\theta \in \mathbb{R}^{\ell}$ and so
$$u^*(t,x)\in H_u(t,x,f_u^{u^*}(t,x)^T\nabla p(t,x)) \quad \mbox{\rm a.e. } (t,x)\in Q_T \, .$$

We get that
$$L(t,x,u_{\varepsilon _n}(t,x))=u_{\varepsilon _n}(t,x)\cdot (f_u^{u_{\varepsilon _n}}(t,x)^T\nabla p_{\varepsilon _n}(t,x)-\varphi (x)(u_{\varepsilon _n}-u^*)(t,x)]$$
$$-H(t,x,f_u^{u_{\varepsilon _n}}(t,x)^T\nabla p_{\varepsilon }(t,x)-\varphi (x)(u_{\varepsilon }-u^*)(t,x))$$
and passing to the limit 
$$L(t,x,u^*(t,x))=u^*(t,x)\cdot f_u^{u^*}(t,x)^T\nabla p(t,x) -\zeta ^*(t,x)$$
a.e. $(t,x)\in Q_T$.
On the other hand
$$L(t,x,u^*)-u^*(t,x)\cdot f_u^{u^*}(t,x)^T\nabla p(t,x)=-H(t,x,f_u^{u^*}(t,x)^T\nabla p(t,x))$$
a.e. and we may conclude that $\zeta ^*(t,x)=H(t,x,f_u^{u^*}(t,x)^T\nabla p(t,x))$.

We may pass to the limit in the equation satisfied by $p_{\varepsilon }$ and get infer that $p$ is a weak solution to (\ref{eqHJp}) and of
\begin{equation*}
	\left\{ \begin{array}{lll}
\displaystyle
\frac{\partial p}{\partial t}(t,x)& \displaystyle +a_{ij}(t,x)\frac{\partial ^2p}{\partial x_i\partial x_j}(t,x)+H(t,x,f_u^{u^*}(t,x)^T \nabla p(t,x)) & \vspace{2mm} \\
~ & \in (f_u^{u^*}(t,x)u^*(t,x)-f^{u^*}(t,x))\cdot \nabla p(t,x)+G_r(t,x,\rho^*(t,x)),
\quad &(t,x)\in Q_T,\vspace{2mm} \\
		p(T,x)&\in -G^0_r(x,\rho ^*(T,x)), \quad &x\in \mathbb{R}^d.
	\end{array}
	\right.
\end{equation*}
\vspace{3mm}

\noindent
\bf Remark. \rm If $f(t,x,u)=f_0(t,x)u$, then $p$ is a weak solution to
\begin{equation*}
	\left\{ \begin{array}{lll}
		\displaystyle
		\frac{\partial p}{\partial t}(t,x)& \displaystyle +a_{ij}(t,x)\frac{\partial ^2p}{\partial x_i\partial x_j}(t,x)+H(t,x,f_0(t,x)^T \nabla p(t,x)) & \vspace{2mm} \\
		~ & \in G_r(t,x,\rho^*(t,x)), \quad &(t,x)\in Q_T,\vspace{2mm} \\
		p(T,x)&\in -G^0_r(x,\rho ^*(T,x)), \quad &x\in \mathbb{R}^d.
	\end{array}
	\right.
\end{equation*}

\section{Uniqueness of the solution to MFG}

Assume that the Hypotheses in Theorem 3.1 (and Section 3) hold and $f$ satisfies (H6).
\vspace{3mm} 

\noindent
\bf Theorem 4.1. \it Assume that for almost any $(t,x)\in Q_T$ we have that
$r\mapsto G(t,x,r)$ is strictly convex on $[0,+\infty )$.
If $(u_1^*,\rho _1,p_1), \ (u_2^*, \rho _2,p_2)$ satisfy
$u_1^*, u_2^*\in {\cal U}$ and $(\rho _k,p_k)$ ($k\in \{ 1,2\} $) is a weak solution to 
\begin{equation}\label{mfg1}
\left\{ \begin{array}{ll}
\displaystyle
\frac{\partial \rho }{\partial t}(t,x)=-\nabla \cdot (f_0(t,x)u_k(t,x)\rho(t,x))+\frac{\partial ^2(a_{ij}\rho )}{\partial x_i\partial x_j}(t,x) , \ & (t,x)\in Q_T, \vspace{2mm} \\
\displaystyle
\frac{\partial p}{\partial t}(t,x)+a_{ij}(t,x)\frac{\partial ^2p}{\partial x_i\partial x_j}(t,x)+H(t,x,f_0(t,x)^T\nabla p(t,x))=\eta _k(t,x) , \ & (t,x)\in Q_T, \vspace{2mm} \\
\rho (0,x)=\rho _0(x), \ p(T,x)=-\eta _{0k}(x), & x\in \mathbb{R}^d ,  
\end{array}
\right.
\end{equation}
where 

$u_k^*(t,x)\in H_q(t,x,f_0(t,x)^T\nabla p_k(t,x))$, a.e. $(t,x)\in Q_T$, 

$\eta_k(t,x)\in G_r(t,x,\rho _k(t,x))$, a.e. $(t,x)\in Q_T$ and 

$\eta _{0k}(x)\in G^0_r(x, \rho_k(T,x))$, a.e. $x\in \mathbb{R}^d$,

\noindent
then $\rho _1\equiv \rho _2$.
\vspace{2mm}

\it Proof. \rm By (\ref{mfg1}) we have that $\rho _1-\rho _2$ is the weak solution to 
\begin{equation}\label{mfg2}
\left\{ \begin{array}{ll}
\displaystyle
\frac{\partial (\rho_1-\rho_2) }{\partial t}(t,x)=-\nabla \cdot (f_0(t,x)(u_1^*\rho_1-u_2^*\rho_2)(t,x)) & \vspace{2mm} \\
\hspace{3cm} \displaystyle +\frac{\partial ^2(a_{ij}(\rho_1-\rho_2)}{\partial x_i\partial x_j}(t,x) , \ & (t,x)\in Q_T, \vspace{2mm} \\
(\rho _1-\rho_2)(0,x)=0, & x\in \mathbb{R}^d   
\end{array}
\right.
\end{equation}
and $p_1-p_2$ is the weak solution to
\begin{equation}\label{mfg3}
\left\{ \begin{array}{ll}
\displaystyle
\frac{\partial (p_1-p_2)}{\partial t}(t,x)+a_{ij}(t,x)\frac{\partial ^2(p_1-p_2)}{\partial x_i\partial x_j}(t,x) & \vspace{2mm} \\
\hspace{3cm} +H(t,x,f_0(t,x)^T\nabla p_1(t,x))-H(t,x,f_0(t,x)^T\nabla p_2(t,x)) & \vspace{2mm} \\
\hspace{3cm} =(\eta _1-\eta_2)(t,x) , \ & (t,x)\in Q_T, \vspace{2mm} \\
(p_1-p_2)(T,x)=-(\eta _{01}-\eta _{02})(x), & x\in \mathbb{R}^d .
\end{array}
\right.
\end{equation}
If we multiply (\ref{mfg2}) by $p_1-p_2$, integrate over $Q_T$ and use (\ref{mfg3}) we obtain that
\begin{equation}\label{mfg4}
\begin{array}{ll}
0=&\displaystyle \int_{\mathbb{R}^d}(\eta_{01}(x)-\eta_{02}(x))(\rho_1(T,x)-\rho_2(T,x))dx  \vspace{2mm} \\
& +\displaystyle \int_{Q_T}(\eta_1(t,x)-\eta_2(t,x))(\rho_1(t,x)-\rho_2(t,x))dx \, dt \vspace{2mm} \\
& +\displaystyle
\int_{Q_T}(u_1^*(t,x)\rho_1(t,x)-u_2^*(t,x)\rho_2(t,x))\cdot f_0(t,x)^T\nabla (p_1-p_2)(t,x)dx \, dt \vspace{2mm} \\
& -\displaystyle
\int_{Q_T}(H(t,x,f_0(t,x)^T\nabla p_1(t,x))-H(t,x,f_0(t,x)^T\nabla p_2(t,x))(\rho _1-\rho_2)(t,x)dx \, dt \, .
\end{array}
\end{equation}
The convexity of $q\mapsto H(t,x,q)$ and $r\mapsto G^0(x,r)$ and the assumptions on $u_k^*$ and $\eta _{0k}$ imply that
$$u_2^*(t,x)\cdot f_0(t,x)^T\nabla (p_1-p_2)(t,x)\rho _2(t,x)\leq (H(t,x,f_0(t,x)^T\nabla p_1(t,x))-H(t,x,f_0(t,x)^T\nabla p_2(t,x))\rho_2(t,x),$$
$$u_1^*(t,x)\cdot f_0(t,x)^T\nabla (p_2-p_1)(t,x)\rho _1(t,x)\leq (H(t,x,f_0(t,x)^T\nabla p_2(t,x))-H(t,x,f_0(t,x)^T\nabla p_1(t,x))\rho_1(t,x),$$
a.e. $(t,x)\in Q_T$,
$$0\leq (\eta_{01}(x)-\eta _{02}(x))(\rho_1(T,x)-\rho_2(T,x)) \quad \mbox{\rm a.e. } x\in \mathbb{R}^d .$$
By (\ref{mfg4}) we get that
$$0\geq \int_{Q_T}(\eta_1(t,x)-\eta_2(t,x))(\rho_1(t,x)-\rho_2(t,x))dx \, dt ,$$
which implies that
$$\rho _1(t,x)=\rho _2(t,x) \quad \mbox{\rm a.e. } (t,x)\in Q_T $$
(because $r\mapsto G(t,x,r)$ is strictly convex). So, we get that $\rho _1=\rho _2=\rho $ in $C([0,T];L^2(\mathbb{R}^d))$
(the uniqueness of $\rho $ from the MFG system).
\vspace{3mm}

\noindent
\bf Remark. \rm If for almost any $(t,x)\in Q_T$ and a.e $x\in \mathbb{R}^d$, respectively, $G$ and $G^0$ are differentiable with respect to $r$ (and so $\eta _1(t,x)=\eta _2(t,x)$
a.e. $(t,x)\in Q_T$ and $\eta _{01}(x)=\eta _{02}(x)$ a.e. $x\in \mathbb{R}^d$), then we get by (\ref{mfg1}) and using the Lipschitz property of $H$ with respect to $q$, that 
$p_1\equiv p_2$. So $p_1=p_2=p$ in $C([0,T];L^2(\mathbb{R}^d))$
(the uniqueness of $p $ from the MFG system).

Moreover, if for almost any $(t,x)\in Q_T$, $H(t,x,\cdot )$ is differentiable, then we get that $u_1^*(t,x)=u_2^*(t,x)$ a.e. $(t,x)\in Q_T$
and so the uniqueness of $u^*$.

\section{Final comments}

Assume that besides the hypotheses in the previous section we get that
$L$ is a continuous function, $F, \, F^0$ are single valued and continuous functions on $\overline{Q}_T\times \mathbb{R}$ and $\mathbb{R}^{d+1}$, respectively, then if $(u^*,\rho ^*)$ is an optimal pair for problem to (P),
there exists a unique weak solution $p$ to
\begin{equation}\label{eqHJG}
\left\{ \begin{array}{ll}
\displaystyle
\frac{\partial p}{\partial t}(t,x)+f^{u^*}(t,x)\cdot \nabla p (t,x)-L(t,x,u^*(t,x)) & \vspace{2mm} \\
\hspace{45mm} \displaystyle +a_{ij}(t,x)\frac{\partial ^2p }{\partial x_i\partial x_j}(t,x)=F(t,x,\rho ^*(t,x)), \quad &(t,x)\in Q_T,\vspace{2mm} \\
p(T,x)+F^0(x,\rho ^*(T,x))=0, & x\in \mathbb{R}^d  
                                     \end{array}
                          \right. 
\end{equation}
(via Lions' existence theorem)
and 
\begin{equation}\label{equ*G}
u^*(t,x)=\mbox{\rm argmin}_{v\in U_0}\{-f(t,x,v)\cdot \nabla p (t,x)+L(t,x,v)\}, \quad \mbox{\rm a.e. } (t,x)\in Q_T \, .
\end{equation}
The proof of this result follows beginning from the fact that
\begin{equation}\label{eqdifG}
0\leq J(u_{\varepsilon }^*,\rho ^{u^*_{\varepsilon }})- J(u^*,\rho ^*) ,
\end{equation}
for any $\varepsilon >0$ and $v\in U_0$, where 
$$u^*_{\varepsilon }(t,x)=\left\{ \begin{array}{ll}
v, \quad &\mbox{\rm if } (t,x)\in B_{\varepsilon }(t_0,x_0) , \vspace{1mm} \\
u^*(t,x), \quad & \mbox{\rm if } (t,x)\in Q_T\setminus B_{\varepsilon }(t_0,x_0) , 
\end{array}
\right. $$
$(t_0,x_0)\in Q_T$ and $\varepsilon >0$ is arbitrary such that
$$B_{\varepsilon }(t_0,x_0) =\{ (t,x)\in \mathbb{R}^{d+1}; \, |(t,x)-(t_0,x_0)| \leq \varepsilon \} \, .$$ 
Note that $u^*_{\varepsilon }$ is a ``spike'' perturbation of $u^*$. 

Dividing (\ref{eqdifG}) by $m(B_{\varepsilon }(t_0,x_0))$ (where $m(A)$ denotes the Lebesgue measure of a set $A$) and arguing as is \cite{anita2021}, \cite{sumin} it follows that for almost any $(t_0,x_0)\in Q_T$:
$$-f^{u^*}(t_0,x_0)\cdot \nabla p (t,x)+L(t_0,x_0,u^*(t_0,x_0))\leq -f(t,x,v)\cdot \nabla p (t,x)+L(t,x,v) \quad \mbox{\rm a.e. } v\in U_0 $$
and so (\ref{equ*G}) holds.

An argument from the first section implies that $p$ is also a weak solution to the next Hamilton-Jacobi equation as well
\begin{equation*}
\left\{ \begin{array}{ll}
\displaystyle
\frac{\partial p}{\partial t}(t,x)+H\left(t,x,f_u^{u^*}(t,x)^T\nabla p(t,x)\right) +a_{ij}(t,x)\frac{\partial ^2p}{\partial x_i\partial x_j}(t,x) & \vspace{2mm} \\
\hspace{25mm} \displaystyle +(f^{u^*}(t,x)-f_u^{u^*}(t,x)u^*(t,x))\cdot \nabla p(t,x)=F(t,x,\rho ^{u^*}(t,x)), \ &(t,x)\in Q_T,\vspace{2mm} \\
p(T,x)+F^0(x,\rho ^*(T,x))=0, & x\in \mathbb{R}^d \, ,
\end{array}
\right.
\end{equation*}
and that
$$u^*(t,x)\in H_q(t,x,f_u^{u^*}(t,x)^T\nabla p (t,x)), \quad (t,x)\in Q_T \, ,$$
where $H(t,x,\cdot )$ is the Hamiltonian of $L(t,x,\cdot )$.
\vspace{3mm}

Let us prove now an uniqueness result for the next MFG system
\vspace{3mm}

\noindent
\bf Theorem 5.1. \it Assume that for almost any $(t,x)\in Q_T$ we have that
$r\mapsto G(t,x,r)$ is strictly convex on $[0,+\infty )$.
If $(u_1,\rho _1,p_1), \ (u_2, \rho _2,p_2)$ satisfy
$u_1, u_2\in {\cal U}$ and $(\rho _k,p_k)$ ($k\in \{ 1,2\} $) is a weak solution to 
\begin{equation}\label{mfg1}
\left\{ \begin{array}{ll}
\displaystyle
\frac{\partial \rho _k}{\partial t}(t,x)=-\nabla \cdot (f^{u_k}(t,x)\rho _k(t,x))+\frac{\partial ^2(a_{ij}\rho _k)}{\partial x_i\partial x_j}(t,x) , \ & (t,x)\in Q_T, \vspace{2mm} \\
\displaystyle
\frac{\partial p_k}{\partial t}(t,x)+f^{u_k}(t,x)\cdot \nabla p _k(t,x)-L(t,x,u_k(t,x)) & \vspace{2mm} \\
\hspace{45mm} \displaystyle +a_{ij}(t,x)\frac{\partial ^2p _k}{\partial x_i\partial x_j}(t,x)=G_r(t,x,\rho _k(t,x)), \quad &(t,x)\in Q_T,\vspace{2mm} \\
\rho _k(0,x)=\rho _0(x), \ p_k(T,x)=-G^0_r(x, \rho _k(T,x)), & x\in \mathbb{R}^d ,  
\end{array}
\right.
\end{equation}
where 
$$u_k(t,x)=\mbox{\rm argmax}_{v\in U_0}\{f(t,x,v)\cdot \nabla p _k(t,x)-L(t,x,v)\}, \quad \mbox{\rm a.e. } (t,x)\in Q_T \, ,$$
then $\rho _1\equiv \rho _2$.
\vspace{2mm}

Proof. \rm If we consider the equation satisfied by $\rho _1-\rho _2$ and the one by $p_1-p_2$, then arguing as in the proof of Theorem 4.1 we get after some calculation that
\begin{equation*}
\begin{array}{ll}
0=&\displaystyle \int_{\mathbb{R}^d}(G^0_r(x,\rho _1(T,x))-G^0_r(x,\rho _2(T,x)))(\rho_1(T,x)-\rho_2(T,x))dx  \vspace{2mm} \\
& +\displaystyle \int_{Q_T}(G_r(t,x, \rho_1(t,x))-G_r(t,x, \rho _2(t,x)))(\rho_1(t,x)-\rho_2(t,x))dx \, dt \vspace{2mm} \\
& +\displaystyle
\int_{Q_T}\Big[ f^{u_1}\rho_2\cdot \nabla p_1+f^{u_2}\rho _1\cdot \nabla p_2-f^{u_1}\rho _1\cdot \nabla p_2 -f^{u_2}\rho _2\cdot \nabla p_1 \vspace{2mm} \\
& +\displaystyle 
L(\cdot ,\cdot ,u_1)\rho_1-L(\cdot ,\cdot, u_1)\rho _2-L(\cdot ,\cdot ,u_2)\rho_1+L(\cdot ,\cdot, u_2)\rho _2\Big] \, dx \ dt .
\end{array}
\end{equation*}
The properties of $u_k$ imply that
$$\rho _2(f^{u_1}\cdot \nabla p_1-L(\cdot ,\cdot ,u_1))\geq \rho _2(f^{u_2}\cdot \nabla p_1-L(\cdot ,\cdot ,u_2)) \, ,$$
$$\rho _1(f^{u_2}\cdot \nabla p_2-L(\cdot ,\cdot ,u_2))\geq \rho _1(f^{u_1}\cdot \nabla p_2-L(\cdot ,\cdot ,u_1)) $$
and by the convexity of $r\mapsto G^0(t,x,r)$ we get that
$$\int_{\mathbb{R}^d}(G^0_r(x,\rho _1(T,x))-G^0_r(x,\rho _2(T,x)))(\rho_1(T,x)-\rho_2(T,x))dx \geq 0 \, .$$
We may infer that
$$\int_{Q_T}(G_r(t,x, \rho_1(t,x))-G_r(t,x, \rho _2(t,x)))(\rho_1(t,x)-\rho_2(t,x))dx \, dt \leq 0$$
and so $\rho _1\equiv \rho _2$ (by the strict convexity of $r\mapsto G(t,x,r)$).
\vspace{3mm}

\noindent
\bf Remark. \rm The last remark in Section 4 clarifies the uniqueness of $p$ and of $u$.
\vspace{3mm}

\noindent
\bf Remark. \rm Theorem 5.1 provides an uniqueness result for the solution to a MFG system for a function $f^u$ which depends nonlinearly on $u$, but under more restrictive 
assumptions on $L, F, F^0$ (which are continuous). 

On the other hand, Theorem 4.1 clarifies the uniqueness of the solution to a more general MFG system (the equation for $p$ is more general) under weaker hypotheses concerning $L, F, F^0$ but for 
function $f^u$ linear with respect to $u$. 
\vspace{3mm}

\noindent
\bf Remark. \rm The existence and uniqueness of the solution to a MFG system for the case when $f\equiv u$ and $a_{ij}=\nu \delta _{ij}$, where $\delta _{ij}$ is the symbol of Kronecker,
has been obtained in \cite{barbu2024} for the more general situation when $\rho _0\in L^1(\mathbb{R}^d)$ (does not necessary belong to $L^2(\mathbb{R}^d)$ as well).
The paper uses some refined regularity arguments and different assumptions on $G$ and $G^0$.

\end{document}